\documentclass[10pt,lettersize,journal]{IEEEtran}
\usepackage{amsmath,amsfonts}
\usepackage{amssymb,amsthm}

\usepackage{algorithm}
\usepackage{algorithmicx}
\usepackage{algpseudocode}
\usepackage{array}

\usepackage{textcomp}
\usepackage{stfloats}
\usepackage{url}
\usepackage{verbatim}
\usepackage{graphicx}
\usepackage{epsfig}
\usepackage{lscape}
\usepackage{subfigure}
\usepackage{epstopdf}

\usepackage{enumerate}
\usepackage{comment,cite,color}
\usepackage{cleveref}
\usepackage{makecell}

\usepackage{bm}
\begin{document}

\renewcommand{\theequation}{\thesection.\arabic{equation}}
\renewcommand{\thefigure}{\arabic{figure}}
\newtheorem{definition}{Definition}[section]
\newtheorem{Proposition}{Proposition}[section]
\newtheorem{theorem}{Theorem}[section]
\newtheorem{corollary}{Corollary}[section]
\newtheorem{example}{Example}[section]
\newtheorem{lemma}{Lemma}[section]
\newtheorem{remark}{Remark}[section]
\newtheorem{notation}{Notation}[section]
\newtheorem{prop}{Proposition}[section]
 \renewcommand{\algorithmicrequire}{\textbf{Input:}}
 \renewcommand{\algorithmicensure}{\textbf{Output:}}
\newcommand{\norm}[1]{\left\|{#1}\right\|_2}
\newcommand{\E}{{\mathbb E}} 
\newcommand{\R}{{\mathbb R}}
\newcommand{\C}{{\mathbb C}}
\newcommand{\T}{{\mathcal{T}}}
\newcommand{\cS}{{\mathcal{S}}} 
\newcommand{\I}{{\mathcal{I}}}
\newcommand{\CN}{{\mathcal{CN}}}
\newcommand{\supp}{{\operatorname{supp}}}
\newcommand{\dist}{{\operatorname{dist}}}
\newcommand{\N}{{\mathcal{N}}}
\newcommand{\eproof}{\hfill\rule{2.2mm}{3.0mm}}
\newcommand{\esubproof}{\hfill$\Box$}
\newcommand{\Proof}{\noindent {\bf Proof.~~}}

\title{GraHTP: A Provable Newton-like Algorithm for Sparse Phase Retrieval}

\author{\IEEEauthorblockN{Licheng Dai\IEEEauthorrefmark{2},
Xiliang Lu\IEEEauthorrefmark{3},
Juntao You\thanks{*Corresponding author: youjuntao@whu.edu.cn.}\IEEEauthorrefmark{1}$^,$\IEEEauthorrefmark{4}$^,$\IEEEauthorrefmark{5}
}

\IEEEauthorblockA{\IEEEauthorrefmark{2}School of Mathematics and Statistics, Wuhan University\\ Wuhan 430072, China \\
\IEEEauthorrefmark{3}School of Mathematics and Statistics, Hubei Center for Applied Mathematics, and Hubei Key Laboratory of Computational Science, Wuhan University\\ Wuhan 430072, China \\ \IEEEauthorrefmark{4}School of Artificial Intelligence, and Hubei Center for Applied Mathematics, Wuhan University, Wuhan 430072, China\\ \IEEEauthorrefmark{5}Institute for Advanced Study, Shenzhen University, Shenzhen 518000, China
}
}

\maketitle

\begin{abstract}
This paper investigates the sparse phase retrieval problem, which aims to recover a sparse signal from a system of quadratic measurements. In this work, we propose a novel non-convex algorithm, termed Gradient Hard Thresholding Pursuit (GraHTP), for sparse phase retrieval with complex sensing vectors. GraHTP is theoretically provable and exhibits high efficiency, achieving a quadratic convergence rate after a finite number of iterations, while maintaining low computational complexity per iteration. Numerical experiments further demonstrate GraHTP's superior performance compared to state-of-the-art algorithms.
\end{abstract}


\section{Introduction}
\IEEEPARstart{P}{hase} retrieval is to recover a signal from the squared modulus of its linear transform, which can be denoted as finding an $n$-dimensional signal $\bm{x}^\dag$ from a system of quadratic equations in the form
\begin{align}\label{eq:problemPR}
    {y}_j = \left|\langle \bm{a}_j, \bm{x}^\dag\rangle\right|^2, j=1,\ldots,m,
\end{align}
where $\{{y}_j\}_{j=1}^{m}\subset \R^m_{+}$ are observed data, $\{\bm{a}_j\}_{j=1}^{m}\subset \C^n $ are given sensing vectors, and $m$ is the number of measurements. The phase retrieval problem arises naturally in fields where direct phase acquisition is difficult or unattainable, such as optics\cite{walther1963question, shechtman2015phase}, X-ray crystallography\cite{harrison1993phase}, quantum mechanics\cite{reichenbach2012philosophic}, quantum information\cite{heinosaari2013quantum}, and others\cite{fienup1987phase,miao2008extending,candes2015phaseconvex,wang2017sparse}.

Solving the nonlinear system~\eqref{eq:problemPR} presents significant challenges. Without additional assumptions on $\bm{x}^\dag$, the system~\eqref{eq:problemPR} may have multiple solutions. Ensuring a unique solution (up to a global phase) requires the so-called ``oversampling'' (i.e. $m>n$) technique, where $m \geq 2n - 1$ for real signals and $m \geq 4n - 4$ for complex signals have been shown to be sufficient with generic sensing vectors \cite{balan2006signal, conca2015algebraic}. Developing practical algorithms for this problem is also highly challenging, which can be traced back to the classical works of Gerchberg-Saxton \cite{gerhberg1972practical} and Fineup\cite{fienup1982phase}. However, earlier approaches lacked rigorous theoretical guarantees. Despite NP-hardness of the problem, a number of practical algorithms that are guaranteed to find true signal (up to global phase) in probabilistic models have been introduced in recent years, which can be categorized into convex and non-convex approaches. Typical convex approaches including PhaseLift\cite{candes2013phaselift,candes2015phaseconvex}, PhaseCut\cite{waldspurger2015phase}, PhaseMax\cite{goldstein2018phasemax} and Flexible convex relaxation\cite{2017A}, enable to recovery signal exactly. However, they can be computationally expensive, especially those that use semi-definite programming (SDP) relaxation and lift the phase retrieval problem to higher dimensional space. Recent studies have introduced several non-convex approaches, including AltminPhase \cite{netrapalli2013phase}, Wirtinger flow \cite{candes2015phase}, truncated amplitude flow (TAF)~\cite{wang2017solving}, Kaczmarz \cite{wei2015solving,tan2019phase}, Riemannian optimization~\cite{cai2024solving},  Gauss-Newton methods \cite{gao2017phaseless, ma2018globally}, among others. These algorithms typically require an initial guess that is sufficiently close to the ground truth to ensure successful recovery, for which the spectral initialization method and its various variants are commonly employed under random Gaussian measurements. Additionally, some approaches have been considered in the context of masked Fourier measurements~\cite{candes2013phaselift,LI2020Phase,li2022sampling}. The number of measurements required by these provable algorithms are $m \sim \mathcal{O}(n \log^a n)$ with $a \geq 0$, which is (nearly) optimal.

Nevertheless, there is significant interest in reducing the necessary number of measurements $m$, especially in high-dimensional applications. This requires leveraging additional information about the unknown signal. In many signal/image processing applications, it is well-established that natural signals or images are often (approximately) sparse in a transformed domain\cite{mallat1999wavelet}. Assuming the $n$-dimensional target signal $\bm{x}^\dag$ is at most $s$-sparse, where $s \ll n$, leads to a sparse phase retrieval problem: recovering $\bm{x}^\dag$ from
\begin{equation}\label{spr}
y_j = \left| \left\langle\bm a_j, \bm{x}^\dag \right\rangle\right|^2,  \, j=1,\dots, m, \quad \text{s.t.} \quad \|\bm{x}^\dag\|_0\leq s,  
\end{equation}
where $\|\bm{x}^\dag\|_0$ denotes the number of nonzero components in $\bm{x}^\dag$. Sparse phase retrieval allows for the recovery of the target signal from an underdetermined system ($m < n$), making it possible to solve phase retrieval problem when only a small number of phaseless measurements are available in practice. It has been shown that $m = \mathcal{O}(s)$ generic measurements are sufficient to determine a unique solution \cite{wang2014phase,akccakaya2015sparse}. However, practical solvers for \eqref{spr} face challenges due to the inherent non-linearity and non-smoothness, especially in underdetermined systems.

\subsection{Related Work and Our Contributions}

The sparse phase retrieval problem \eqref{spr}
has been extensively studied in recent years, and a number of provable practical algorithms have been introduced. For instance, the $\ell_1$-regularized PhaseLift \cite{li2013sparse}, a natural extension of the convex approach to the compressive case, demonstrates that the signal $\bm{x}^\dag$ can be correctly recovered with $\mathcal{O}(s^2\log n)$ random Gaussian measurement. Non-convex algorithms have also garnered lots attention due to their computational efficiency. Non-convex approaches typically involve two stages: an initialization stage followed by a local refinement stage. Such algorithms include 
  SPARTA\cite{wang2017sparse}, CoPRAM\cite{jagatap2019sample}, thresholding/projected
Wirtinger flow \cite{cai2016optimal,soltanolkotabi2019structured},   SAM\cite{cai2022sample}, HTP\cite{cai2022sparse} and others\cite{cai2024a}. Given an initial guess sufficiently close to the target signal, most of these algorithms guarantee at least linear convergence to the ground truth with $\mathcal{O}(s\log (n/s))$ random Gaussian measurements.  Meanwhile, an initial guess sufficiently close to the ground truth can be generated using specific methods, such as spectral methods, which can produce an appropriate initial guess using $m \sim \mathcal{O}(s^2\log n)$ random Gaussian measurements\cite{netrapalli2013phase,wang2017sparse}. For a more detailed discussion, see \cite{cai2023provable}.

Theoretically, algorithms like ThWF, SPARTA and CoPRAM achieve $\varepsilon$-accuracy in $\mathcal{O} ( \log(1 / \varepsilon))$ iterations. Recently proposed methods such as SAM\cite{cai2022sample} and HTP\cite{cai2022sparse} guarantee exact recovery of target signal within a finite number of iterations, specifically at most $\mathcal{O}(\log(s^2 \log n) + \log(\|\bm{x}^\dag\|_2 / x^\dag_{\min}))$ iterations. The HT-Newton (hard-thresholding Newton) method in \cite{cai2024a} establishes a quadratic convergence rate to the ground truth (up to a global sign) after at most $\mathcal{O}(\log (\| \bm {x}^{\dag}\|_2/ {x}^{\dag}_{\min}))$ iterations. Those work demonstrate improved efficiency in solving the sparse phase retrieval problem. However, these theories of finite-step/superlinear convergence are both based either on the case of real-valued $\bm{x}^\dag$ and $\{\bm{a}_j\}_{j=1}^m$, or on assumptions that certain problems on subspaces can be exactly solved \cite{xu2024subspace}. It is known that complex sensing vectors $\{\bm{a}_j\}_{j=1}^m$ are of great interest in the applications\cite{candes2013phaselift}.  In this work, we explore the scenario where the sensing vectors are complex and present a novel, efficient Newton-like algorithm which is guaranteed to achieve superlinear convergence. Our main contributions are four-fold:
\begin{itemize}
    \item We extend the Gradient Hard Thresholding Pursuit (GraHTP) algorithm \cite{yuan2016exact, yuan2018gradient} to minimize the intensity-based loss for sparse phase retrieval with a complex sensing matrix. The proposed algorithm is highly efficient: to find an $\varepsilon$-solution, the number of iterations is at most $\mathcal{O}(\log (\log (1/{\varepsilon})) + \log({\|\bm{x}^\dag }\|_2/{{x}^\dag_{\min}})) $ with a per iteration complexity $\mathcal{O}(mn+L_0 s^2m)$, $L_0\ge 1$ is any fixed constant. For $s\ll n$, the per-iteration complexity is of the same order as that of first-order gradient-type methods.  
    \item The theoretical guarantee of quadratic convergence rate, achieved after at most $\mathcal{O}(\log (\|\bm{x}^\dag\|_2/ x_{\min}^\dag))$ iterations, has been established for GraHTP under certain conditions. 
    \item With the resampling technique, we prove that each local update of GraHTP requires only
$m=\mathcal O(s\log(n/s))$ samples. Combined with the modified spectral initialization
\cite{cai2023provable}, this gives a total sample complexity
$\mathcal O(\bar s\,s\log n+K_\varepsilon s\log(n/s))$, where
$\bar s=\|\bm x^\dagger\|_2^2/\|\bm x^\dagger\|_\infty^2$ and
$K_\varepsilon=\mathcal O(\log\log(1/\varepsilon)+
\log(\|\bm x^\dagger\|_2/x^\dagger_{\min}))$.
    \item The empirical advantages of the proposed algorithm have been verified against state-of-the-art algorithms. Numerical experimental results illustrate that GraHTP achieves the superior computational efficiency and recoverability in the complex test problems.
\end{itemize}

\subsection{Notation}
For an index set $\cS \subseteq \{1,2,\cdots,n\}$, $|\cS|$ denotes the number of elements in set $\cS$. For a vector $\bm{z} \in \R^n$, $\norm{\bm{z}} $ denotes the Euclidean norm. $\bm{z}_\cS$ means the sub-vector indexed by $\cS$.  For a matrix $\bm{A} \in \R^{m\times n}$, $\norm{\bm{A}}$ denotes the spectral norm, $\bm{A}_\cS$ represents retaining only the columns of the matrix indexed by $\cS$, and $\bm{A}_{\cS,\T}$ represents retaining the columns and rows of the matrix indexed by $\cS$ and $\T$ respectively. 
$\bm{I}_{|\T|} $ denotes the $|\T|$ dimensional identity matrix.   Hard thresholding operator $\mathcal{H}_s: \R^n \rightarrow \R^n$ represents retaining the maximum $s$ components in the magnitude of a vector in $\R^n$, and setting the other components to zero. $x^\dag_{\min}$ represents the nonzero component with the smallest absolute value of the vector $\bm{x}^\dag$.

\section{Algorithms}
In this section, we describe our proposed algorithm in detail. The proposed algorithm is based on the general framework of GraHTP\cite{yuan2016exact,yuan2018gradient}, which enjoys finite-step convergence in the case of compressive sensing\cite{foucart2011hard}. And similar to most of the existing non-convex
sparse phase retrieval algorithms, the proposed algorithm consists of two stages, namely, the initialization stage and the iterative refinement stage. In this work, we focus on the iterative refinement stage, and the initialization stage can be done by an off-the-shelf algorithm such as the spectral method or modified spectral method\cite{cai2023provable}.

\subsection{The Proposed Algorithm}\label{the proposed algorithm}                                                                                                             
 
In practice, the unknown signal may belong to either $\mathbb{R}^n$ or $\mathbb{C}^n$. For simplicity, we consider the real-valued case $\bm{x}^\dag \in \mathbb{R}^n$. The proposed algorithm can be extended to complex-valued signals by separating the real and imaginary parts \cite{gao2017phaseless}.
The squared error loss associated with \eqref{spr}, commonly referred to as the intensity-based loss in the context of the phase retrieval problem, naturally leads to the following optimization problem:
\begin{align}\label{problem:minspr}
	\min _{\bm{z} \in \R^n} f(\bm{z}) \text {, s.t. }\|\bm{z}\|_0 \leq s,
\end{align}
where 
\begin{align}\label{intensity-based}
      {f}\left(\bm{z}\right) :=\frac{1}{4 m} \sum_{j=1}^m\left(|\langle \bm{a}_{j}, \bm{z}\rangle |^2-{y}_j\right)^2.
  \end{align} 
To address the problem formulated in~\eqref{problem:minspr}, we consider employing the GraHTP method, as introduced in \cite{yuan2016exact,yuan2018gradient}. GraHTP extends the HTP algorithm from compressed sensing to tackle the broader framework of sparsity-constrained convex optimization problems. The approach begins with a projected gradient descent (PGD), followed by subspace selection, and then resolves a subspace optimization problem to refine the solution. Specifically, for a given current estimate $\bm{z}^{k}$ at the $k$-th iteration, the one-step update of $\bm{z}^{k}$ consists of the following three sub-steps:
\begin{itemize}
    \item[1)] Compute the PGD update as 
    \begin{align}\label{u=Hs}
      \bm{u}^{k} = \mathcal{H}_{s}\big(\bm{z}^k-\mu^k\nabla{f}(\bm{z}^k)\big), ~ \text{step size}~ \mu^k>0.
    \end{align}  
    \item[2)] Estimate guess of the support
    $$\cS_{k+1} = \supp(\bm{u}^{k}).$$ 
    \item[3)] Compute update $\bm{z}^{k+1}=\hat{\bm{z}}^{k+1}$ where
     \begin{align}\label{S_k minimize}	
  \hat{\bm{z}}^{k+1} \leftarrow \underset{\supp(\bm{z}) \subseteq 	\cS_{k+1}}{\arg \min } f(\bm{z}).
\end{align}
\end{itemize}
However, applying GraHTP directly to our problem is not feasible both algorithmically and theoretically. Firstly, the exact solution for the optimization problem \eqref{S_k minimize} associated with the loss function in \eqref{intensity-based} is nontrivial to obtain. Moreover, ensuring convergence in the case of non-convex $f(\bm{z})$ and the lack of an exact solution for \eqref{S_k minimize} pose significant challenges. 

For the first issue, it is important to note that we are minimizing $f(\bm{z})$ on restricted supporting set $\cS_{k+1}$, which is small in scale compared to the dimension $n$. Therefore, we propose solving \eqref{S_k minimize} approximately with $L$ steps of Gauss-Newton iteration \cite{bjorck1996numerical,fletcher2000practical}. Let $\{\bm{z}^{k,l}\}_{l=0}^L$ be the sequence generated by the iteration with an initial guess $\bm{z}^{k,0}$. For example, one can choose the initial guess $\bm{z}^{k,0}=\bm{u}^{k}$. Next we perform Gauss-Newton update. For the ease of notation, we let $F_j(\bm{z}):=\frac{1}{2\sqrt{m}}\big(|\langle \bm{a}_{j}, \bm{z}\rangle |^2-y_j\big)$ as the $j$-th component of $\bm{F}(\bm{z}) : \R^n \rightarrow \R^m$, and the loss function $f(\bm{z})$ can be written in the form of:
\begin{align}\label{sparse minimize}
f(\bm{z})=  \sum_{j=1}^{m}{F}_j(\bm{z})^2. 
\end{align}
Performing a first-order Taylor expansion of $\bm{F}(\bm{z})$ at $\bm{z}^{k,l}$:
\begin{align}\label{linearlize}
			\bm{F}(\bm{z}) &\approx  \bm{F}(\bm{z}^{k,l})+\bm{J}(\bm{z}^{k,l})(\bm{z}-\bm{z}^{k,l}) ,
		\end{align}
where the j-th row of matrix $\bm{J}(\bm{z}^{k, l}) \in \R^{m\times n}$ is $\frac{1}{\sqrt{m}}(\bm{a}_{j R} \bm{a}_{j R}^{\top} \bm{z}^{k, l} +\bm{a}_{j I} \bm{a}_{j I}^{\top}\bm{z}^{k, l})^{\top}$ and $\bm{a}_{j R}$, $\bm{a}_{j I}$ are the real and imaginary part of $\bm{a}_j$ respectively. Combining with \eqref{sparse minimize}, the problem in \eqref{S_k minimize} can be approximated by
\begin{equation}\label{approximation}
\underset{\supp(\bm{z}) \subseteq 	\cS_{k+1}}{ \min }  \norm{\bm{J}(\bm{z}^{k,l})(\bm{z}-\bm{z}^{k,l})+\bm{F}(\bm{z}^{k,l})}^2 .
		\end{equation}
Denote $ \bm{z}^{k,l+1}$ as the solution to \eqref{approximation},
which satisfies 
$$\bm{z}^{k,l+1}_{\cS_{k+1}^c} = \bm{0}$$ 
and
\begin{align}\label{iteration3}
 &\bm{J}_{\cS_{k+1}}(\bm{z}^{k,l})^{\top} \bm{J}_{\cS_{k+1}}(\bm{z}^{k,l}) \big(  \bm{z}^{k,l}_{\cS_{k+1}} - \bm{z}^{k,l+1}_{\cS_{k+1}} \big)\cr
			=& \bm{J}_{\cS_{k+1}}(\bm{z}^{k,l})^{\top} \bm{F}(\bm{z}^{k,l}).  
\end{align}  
Finally we choose the update $\bm{z}^{k+1}$ as $\bm{z}^{k,L}$. The proposed algorithm, coined GraHTP for sparse phase retrieval, is summarized in \Cref{algorithm}. Also, the proposed algorithm can be extended to the case of complex $\bm{x}^\dag$ in numerical experiments.\par
For the proposed GraHTP, computing $\bm{z}^{k,0}$ and $\bm{z}^{k+1}$ costs $\mathcal{O}(mn)$ and $\mathcal{O}(s^2m)$ flops, respectively. The total computational cost for updating $\bm z$ is $\mathcal{O}(mn + Ls^2m)$. 
As long as $s\ll n$ or $s\le \sqrt{n}$, the per-iteration complexity of the proposed algorithm is comparable to that of first-order gradient-type methods like thresholding/projected Wirtinger flow \cite{cai2016optimal,soltanolkotabi2019structured}, which also costs $\mathcal{O}(mn)$ flops per-iteration. However, our proposed algorithm guarantees super-linear convergence under certain conditions, as established in the next section.

  \begin{algorithm}
	\caption{\textbf{Gra}dient \textbf{H}ard \textbf{T}hresholding \textbf{P}ursuit for Sparse Phase Retrieval}	\label{algorithm}
 \begin{algorithmic}
     \Require 
  Data $\{\bm{a}_j, {y}_{j}\}_{j=1}^{m}$, the sparsity level $s$, the maximum number $K$ and $L$ of iterations allowed, step size $\mu^k$.
     \begin{enumerate}
         \item[1:] Initialization: Let the initial value $\bm{z}^{0}$ be generated by a given method, e.g., (modified) spectral method\cite{jagatap2019sample,cai2023provable}.
         \item[2:] \textbf{for} $k=0,1,\ldots,K-1$ \textbf{do}
         \item[3:]  $  \bm{u}^{k}=\mathcal{H}_{s} \big (\bm{z}^k - \mu^k\nabla {f}(\bm{z}^k)\big)$
         \item[4:] $\cS_{k+1}=\supp(\bm{u}^{k})$
         \item[5:] Obtain $\bm{z}^{k+1}$ by solving
         \begin{align*}
        \underset{\supp(\bm{z}) \subseteq 	\cS_{k+1}}{ \min }  \sum_{j=1}^{m}{F}_j(\bm{z})^2
         \end{align*}     
         via $L$ steps of Gauss-Newton iteration: starting from  $\bm{z}^{k, 0} = \bm{u}^{k}$ where $\supp(\bm{z}^{k, 0})=\cS_{k+1}$ ,
         
         \quad\textbf{for} $l = 0,\cdots,L-1$ \textbf{do}
             \begin{align*}
            \bm{z}^{k, l+1} = \underset{\supp(\bm{z}) \subseteq 	\cS_{k+1}}{\arg \min } \norm{ \bm{J}(\bm{z}^{k, l})(\bm{z}-\bm{z}^{k, l})+\bm{F}(\bm{z}^{k, l}) }^2 
        \end{align*}       
        \quad\textbf{end for}
        
        Set $\bm{z}^{k+1} = \bm{z}^{k, L}$

 \item[6:] \textbf{end for}
     \end{enumerate}
	\Ensure
 $\bm{z}_{output}=\bm{z}^{K}$.
 \end{algorithmic}	
\end{algorithm}

\subsection{Local Convergence}
In this subsection, we present the theoretical results of the proposed GraHTP for sparse phase retrieval as summarized in \Cref{algorithm}, and we take the number of inner iterations $L$ to be a fixed constant $L_0$. The distance between  $\bm{x}^\dag$ and $\bm{z}$ is define as $\mathrm{dist}(\bm{x}^\dag,\bm{z}) := \min\{\|\bm{x}^\dag-\bm{z}\|_2,\|\bm{x}^\dag+\bm{z}\|_2\}$ as $\bm{x}^\dag$ and $-\bm{x}^\dag$ are equivalent solutions. If the current guess $\bm{z}^k$ fall within a nearby local neighborhood of $\bm{x}^\dag$ or $-\bm{x}^\dag$, the following \Cref{local convergence} demonstrates the one-step contraction property of GraHTP in \Cref{algorithm}, referring to the number of Gauss-Newton iterations. 
For ease of presentation, we denote the basin of attraction as
\begin{equation*}
\mathcal{E}(\gamma):=\{\bm{z}\in\mathbb{R}^n~|~\dist\left(\bm{z}, \bm{x}^\dag\right) \leq \gamma \|\bm{x}^\dag \|_2,  \|\bm{z} \|_0 \leq s\}.
\end{equation*}
for $\gamma \ge 0$.

\begin{theorem}
\label{local convergence}
Let $\bm{x}^\dag \in \R^n$ be any s-sparse signal. Consider $m$ i.i.d. random vectors $\{\bm a_j\}_{j=1}^{m}$ with $\bm{a}_j\sim \CN(\bm{0},\bm{I})$ at every update, including both the outer hard-thresholding step and the inner Gauss--Newton refinement steps, and noiseless measurements $y_j=|\langle\bm{a}_j,\bm{x}^\dag\rangle|^2$ satisfying system \Cref{spr}. Then for some small constant $\gamma \in (0,1)$ and  any $L_0 \geq 1$, there exist universal constants $\mu_1, \mu_2, C_1, C_2$, $\alpha\in (0,1)$ and $\beta$ such that: For any fixed $\bm{z}^k \in \mathcal{E}(\gamma)$, with probability at least $1-C_1 {m}^{-1}$, the $\bm z^{k+1}$ produced by \Cref{algorithm} 
satisfies
\begin{align*}
\dist\left(\bm{z}^{k+1}, \bm{x}^\dag\right)
\begin{cases}
\le \alpha \cdot \dist\left(\bm{z}^k, \bm{x}^\dag\right),
&  
\cr
\le \beta \cdot \dist^2\left(\bm{z}^k, \bm{x}^\dag\right), & \text{if} ~\bm{z}^{k}\in\mathcal{E}(\frac{x^\dag_{\min}}{\norm{\bm{x}^\dag }}),
\end{cases}
\end{align*}
if step size $\mu^k \in \big(\frac{\mu_1}{\norm{\bm{x}^\dag}^2}, \frac{\mu_2}{\norm{\bm{x}^\dag}^2}\big)$ and	$m \geq  C_2 s \log (n / s)$.
\end{theorem}

\begin{proof}
The proof of this theorem is deferred to \Cref{proof3}.
\end{proof}

Therefore, $L_0$ can be chosen as any fixed small constant, such as $L_0=1$ or $L_0=2$. 
It can be found that \Cref{algorithm} exhibits two phases of convergence: initially, in the first phase, a linear rate of convergence is achieved when the current estimate is within an $\mathcal{O}(1) \|\bm{x}^\dag \|_2 $ neighborhood of $\bm{x}^\dag$; subsequently, in the second phase, the algorithm attains a quadratic convergence rate when the current estimate is within an $x^\dag_{\min} $ neighborhood of $\pm\bm{x}^\dag$. 
The linear-to-superlinear transition occurs after
$\dist(\bm z^k,\bm x^\dagger)\leq x^\dagger_{\min}$, and therefore may require
$\mathcal O(\log(\|\bm x^\dagger\|_2/x^\dagger_{\min}))$ linear iterations for ill-conditioned signals.
The resampling assumption can be implemented by collecting all measurements once and partitioning them into independent batches; a fully fixed-data theory is left for future work.

 \begin{remark}
Let $\widehat\nu:=m^{-1}\sum_{j=1}^m y_j$. Under Gaussian sensing, for any $0<\delta<1$, we have
\[
\mathbb P\left\{
(1-\delta)\|\bm x^\dagger\|_2^2
\leq
\widehat\nu
\leq
(1+\delta)\|\bm x^\dagger\|_2^2
\right\}
\geq
1-2\exp(-cm\delta^2).
\]
Hence $\mu^k=\bar\mu/\widehat\nu$ satisfies the stepsize condition in Theorem~II.1 with high probability provided $m\ge C\delta^{-2}\log n$ for some $C>0$.
\end{remark}
 Now we consider the noisy case where the measurements are given by $\bm{y}^{(\eta)}=\bm{y}+\bm{\eta}$, where $\bm{\eta} \in \R^m$ is a noise vector independent of $\{\bm{a}_j \}_{j=1}^m$.
\begin{theorem}\label{noisy case}
Let $\bm{x}^\dag \in \R^n$ be any s-sparse signal. Consider $m$ i.i.d. random vectors $\{\bm a_j\}_{j=1}^{m} $ with $\bm{a}_j\sim \CN(\bm{0},\bm{I})$ at every update, including both the outer hard-thresholding step and the inner Gauss--Newton refinement steps, and noisy measurements $y_j^{(\eta)}=|\langle\bm{a}_j, \bm{x}^\dag\rangle|^2 + \eta_j$. Then for some small constant $\gamma \in (0,1)$, $q_0$ and any $L_0 \geq 1$, there exist universal constants  $\mu_1, \mu_2,  C_{\eta},  C_3, C_4$, $ p \in (0,1)$ and $q \leq q_0$ such that: For any fixed $\bm{z}^k \in \mathcal{E}(\gamma)$, with probability at least $1-C_3 {m}^{-1}$, the $\bm{z}^{k+1}$ produced by \Cref{algorithm} satisfies
\begin{align*}
\dist\left(\bm{z}^{k+1}, \bm{x}^\dag\right) 
\le p \cdot \dist\left(\bm{z}^k, \bm{x}^\dag\right) + q \cdot \norm{\bm{\eta}}, 
 \end{align*}
 provided $\|\bm{\eta}\|_2 \leq C_{\eta} \|\bm{x}^\dag\|_2$, step size $\mu^k \in \big(\frac{\mu_1}{\norm{\bm{x}^\dag}^2}, \frac{\mu_2}{\norm{\bm{x}^\dag}^2}\big)$ and $m \geq  C_4 s \log (n / s)$.
\end{theorem}
\begin{proof}
The proof of this theorem is deferred to \Cref{proof5}.
\end{proof}

The local result above implies that if there exists an initial guess satisfying $\bm{x}^0 \in \mathcal{E}(\gamma)$ provided $m$ at least $\mathcal{O} (s\log(n/s))$, then with overwhelming probability we have
 \begin{align*}
\dist\left(\bm{z}^k, \bm{x}^\dag\right) 
\le \gamma p^k \cdot \|\bm{x}^\dag\|_2 + q \frac{1-p^k}{1-p}\cdot \norm{\bm{\eta}}  , \quad
 \forall k \geq 0 .
 \end{align*}
 
 \subsection{Initialization and Global Convergence}\label{initialization}
The desired initial guess can be produced by algorithms such as spectral method or modified spectral method for random Gaussian measurements. 
The spectral initialization constructs a matrix $\bm{Y} := \frac{1}{m} \sum_{i=1}^m y_i \bm{a}_i \bm{a}_i^*$ or its variants \cite{netrapalli2013phase,candes2015phase,chen2015solving}, whose leading eigenvector is a good approximation to $\pm \bm{x}^\dag$ if the number of Gaussian measurements is at least $\mathcal{O} (n)$. For sparse phase retrieval, the number of Gaussian measurements can be reduced by first estimating the support of $\bm{x}^\dag$ as $\mathcal{S}_0$, e.g., top-$s$ entries in the diagonal elements of $\bm{Y}$ given by $\left\{\frac{1}{m}\sum_{i=1}^{m}y_i |a_{ij}|^2\right\}_{j=1}^{n}$, and the support vector of the initial guess is estimated as the principal eigenvector of the real part of $\frac{1}{m}\sum_{i=1}^my_i\bm{a}_{i,\mathcal{S}_0}\bm{a}_{i,\mathcal{S}_0}^*$. The initial guess $\bm{z}^0$ generated by this setting can be sufficiently close to the ground truth: For any $\gamma \in (0,1)$, with probability at least $1-  8{m}^{-1}$, 
\begin{align}\label{x^0}
\dist\left(\bm{z}^{0}, \bm{x}^\dag\right) \leq \gamma \norm{\bm{x}^\dag },    
\end{align}
provided ${m} \sim \mathcal{O} (s^2 \log n)$, see~Theorem IV.1 in \cite{jagatap2019sample}.

Moreover, the number of measurements required for the initialization stage can be further reduced by a modified version of spectral initialization. For more details, refer to \cite{cai2023provable}, see more discussion in Remark \ref{remark:3}.

\begin{theorem}
\label{global}
Let $\bm{x}^\dag \in \R^n$ be any s-sparse signal. 
For constants $L_0\geq 1$ and $K_{\max} \geq 1$, let the initial guess $\bm{z}^0$ generated by the spectral method described above and the sequence $\{\bm{z}^k\}_{k=1}^{K_{\max}}$  generated by \Cref{algorithm} with independent $\{\bm a_j, y_j\}_{j=1}^{m}$ at every iteration, where $\bm{a}_j\sim \CN(\bm{0},\bm{I})$ and $y_j=|\langle\bm{a}_j, \bm{x}^\dag\rangle|^2$. Then there exist positive constants $\mu_1, \mu_2, C_5, C_6, C_7, C_8$, $\alpha\in (0,1)$, and $\beta$ such that:  with probability at least $1- C_5 {m}^{-1}$, 
\begin{align*}
\dist\left(\bm{z}^{k+1}, \bm{x}^\dag\right)
\begin{cases}
\le \alpha \cdot \dist\left(\bm{z}^k, \bm{x}^\dag\right), & ~ 0\leq k <  K_0
\cr
\le \beta \cdot \dist^2\left(\bm{z}^k, \bm{x}^\dag\right), &  ~ K_0 \leq k < K_{\max}
\end{cases}
\end{align*}
where $K_0 \leq C_6\log\big({\|\bm{x}^\dag\|_2}/{x^\dag_{\min}}\big)+ C_7$, provided step size $\mu^k \in \big(\frac{\mu_1}{\norm{\bm{x}^\dag}^2}, \frac{\mu_2}{\norm{\bm{x}^\dag}^2}\big)$ and	$m \geq  C_8s\log (n/s)$.
\end{theorem}

\begin{proof}
 The proof of this theorem is deferred to \Cref{proof4}.
\end{proof}

\begin{remark}
The linear convergence achieved in the first phase ensures that reaching $x^\dag_{\min}$-closeness to $\pm \bm{x}^\dag$ typically requires at most $\mathcal{O}\big(\log \big(\|\bm{x}^\dag \|_2/ x^\dag_{\min}\big)\big)$ iterations. Consequently, the proposed algorithm realizes a quadratic convergence rate after at most $\mathcal{O}\big(\log \big(\|\bm{x}^\dag \|_2/ x^\dag_{\min}\big)\big)$ iterations. Moreover, to achieve an $\varepsilon$-solution, the required number of iterations $K_{\varepsilon}$ is at most $\mathcal{O}\big(\log (\log (1/{\varepsilon})) + \log({\|\bm{x}^\dag \|_2}/{{x}^{\dag}_{\min}})  \big)$.
Since initialization required $ \mathcal{O}(s^2\log n)$ random Gaussian measurements, according to \Cref{global}, the total sampling complexity required in \Cref{algorithm} is 
    $\mathcal{O}\big(s^2\log n +  K_{\varepsilon}s\log(n/s)\big)$. 
\end{remark}

\begin{remark}\label{remark:3}
Theoretically, the sampling complexity of using\cite{cai2023provable} in the initialization stage can be reduced to $\mathcal{O}(\bar{s}s\log n)$, where $\bar s :=\|\bm x^\dagger\|_2^2/ \|\bm x^\dagger\|_\infty^2$ is the stable sparsity of the underlying signal. This leads to a lower global sampling complexity $\mathcal{O}\big(s\log n +  K_{\varepsilon} s\log(n/s)\big)$ when the signal contains $\mathcal{O}(1)$ significant components.    
\end{remark}

\begin{remark}\label{complex case}
 In the case of complex-valued signals, we followed the programming described in \cite{gao2017phaseless} and extended it to sparse signals. 
 The numerical experiments show that our proposed algorithm performs well
in this case. However, while the sign of real-valued signal is discrete, in the complex-
valued case it varies continuously, which leads to technical limitations in theoretical
guarantees. Addressing this problem will be an important topic for our future research.   
\end{remark}

\section{Numerical Experiments}
In numerical simulation process, the true signal $\bm{x}^\dag$ is set to have $s$ nonzero entries. In the first and second subsections, the sensing vectors $\left\{\bm{a}_j\right\}_{j=1}^{m}$ follow the Gaussian random distribution, i.e., $\bm{a}_j \sim \CN(0, \bm{I})$. In the last subsection, we use partial discrete Fourier transform matrix as the sensing matrix $\bm{A} = [\bm{a}_1, \bm{a}_2, \cdots, \bm{a}_m]^\top \in \C^{m\times n}$. All ground-truth signals were normalized such that $\|\bm x^\dagger\|_2 = 1$. The support set of $\bm{x}^\dag$ is uniformly and randomly extracted from all $s$-subsets of set $\{1,2\cdots,n\}$ and the values of nonzero terms are independently and randomly generated from the standard Gaussian distribution $\N(0, \bm{I})$. $\left\{y_j\right\}_{j=1}^{m}$ are samples without noise, where $y_j=|\langle\bm{a}_j, \bm{x}^\dag\rangle|^2, \quad j=1,2, \cdots, m$. The observation data with noise are defined by the following equation: 
\[
y_j^{(\eta)}=y_j+\sigma\eta_j, \quad j = 1,\cdots,m,
\]
where noise $\left\{\eta_j\right\}_{j=1}^{m}$ follows the standard Gaussian random distribution, and we use $\sigma > 0$ to determine the noise level .

\begin{figure}[ht!]
     \centering
\includegraphics[width=0.21\textwidth,trim=3pt 3pt 30pt 5pt,clip]{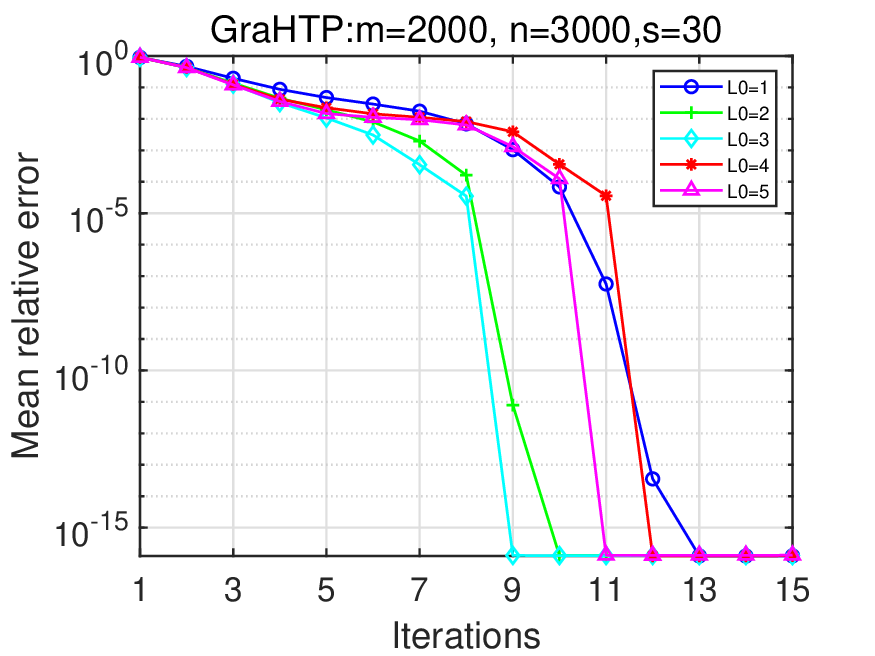}  \includegraphics[width=0.21\textwidth,trim=3pt 3pt 30pt 5pt,clip]{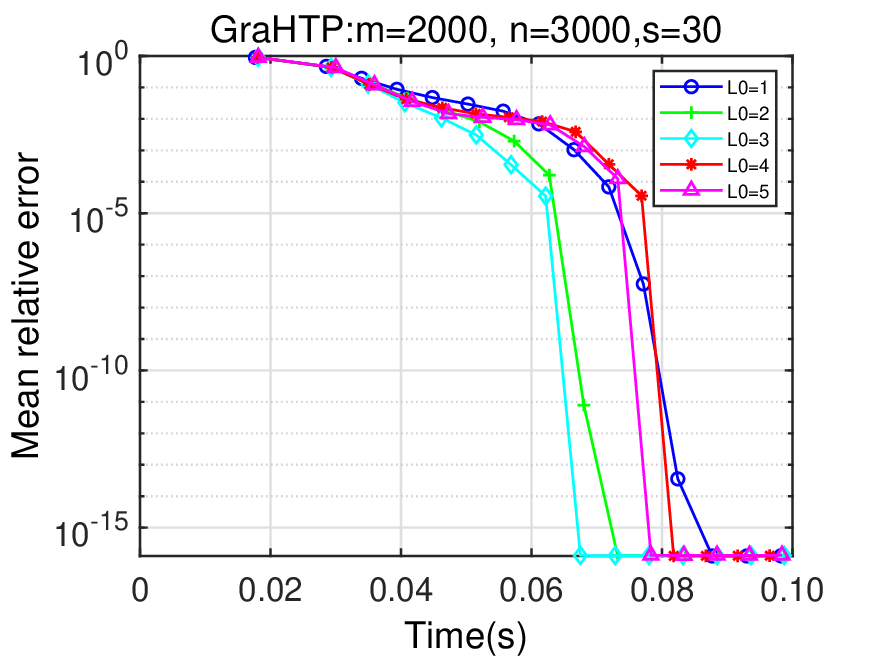} 
     \caption{Relative error versus number of iterations and relative error versus of running time of GraHTP with complex-valued sensing matrix, the results were obtained by averaging 100 independent trials.}
     \label{fig:complexsignal}
   \end{figure}

We will compare our algorithm GraHTP with other popular algorithms such as CoPRAM\cite{jagatap2019sample}, ThWF\cite{cai2016optimal}, SPARTA\cite{wang2017sparse},  HTP\cite{cai2022sparse}, and HT-Newton \cite{cai2024a}. 
We follow the standard practice and use the spectral initialization method to ensure fairness. The numerical experiments are run on a computer with 3.00 GHz Intel Core i9 processor and 64 GB RAM using MATLAB R2023a. In experiments, the parameters of ThWF are set as recommended in \cite{cai2016optimal}. The parameters of SPARTA are set to be $\mu = 1$, $\delta = 0.7$ and $|\I| = \lceil m/6 \rceil$. The step sizes $\mu$ of HTP and HT-Newton are both fixed to $0.95$. For GraHTP, we set the normalized step size to $\bar\mu=0.35$ and fix the
number of inner Gauss--Newton iterations as $L_0=3$. Similar performance was
observed for other small fixed choices of $L_0$ (see \Cref{fig:complexsignal}). The relative error between the true signal $\bm{x}^\dag$ and the estimated signal $\hat{\bm{x}}$ is defined as
\begin{equation}
	r\left(\hat{\bm{x}}, \bm{x}^\dag\right)=\frac{\dist\left(\hat{\bm{x}}, \bm{x}^\dag\right)}{\norm{\bm{x}^\dag}} .
\end{equation}
where $\dist\left(\hat{\bm{x}}, \bm{x}^\dag\right) = \min_{\phi\in[0,2\pi)}\norm{ \hat{\bm{x}}-e^{i\phi}\bm{x}^\dag }$. We define that signal recovery is successful when $r\left(\hat{\bm{x}}, \bm{x}^\dag\right) \leq 10^{-6}$. 

\subsection{Real-valued Signal Case}

\begin{figure}[ht!]
\begin{center}
\subfigure[sparsity $s = 20$]{
\includegraphics[width=0.22\textwidth,trim=3pt 3pt 30pt 5pt,clip]{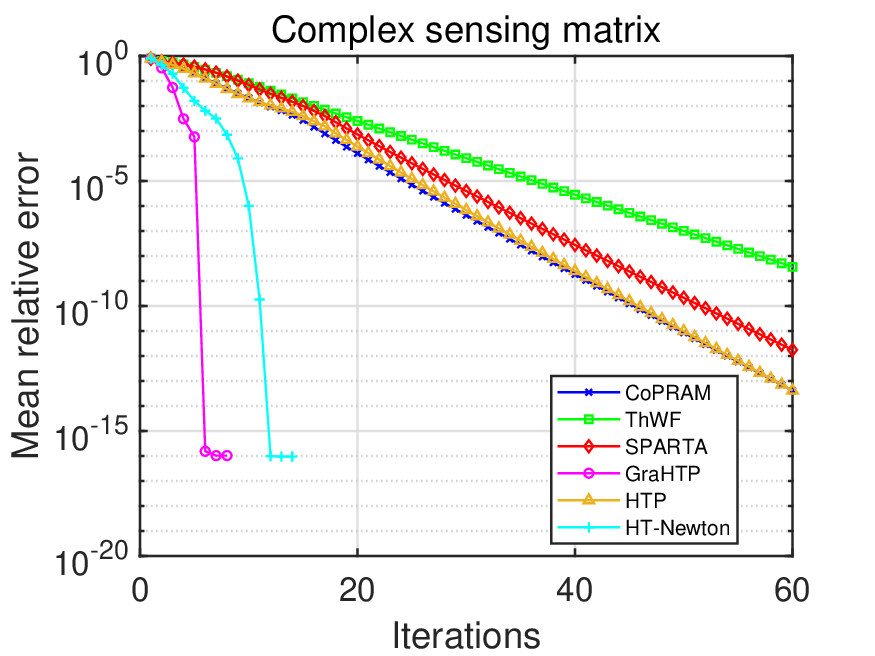}}
\subfigure[sparsity $s = 30$]{
\includegraphics[width=0.22\textwidth,trim=3pt 3pt 30pt 5pt,clip]{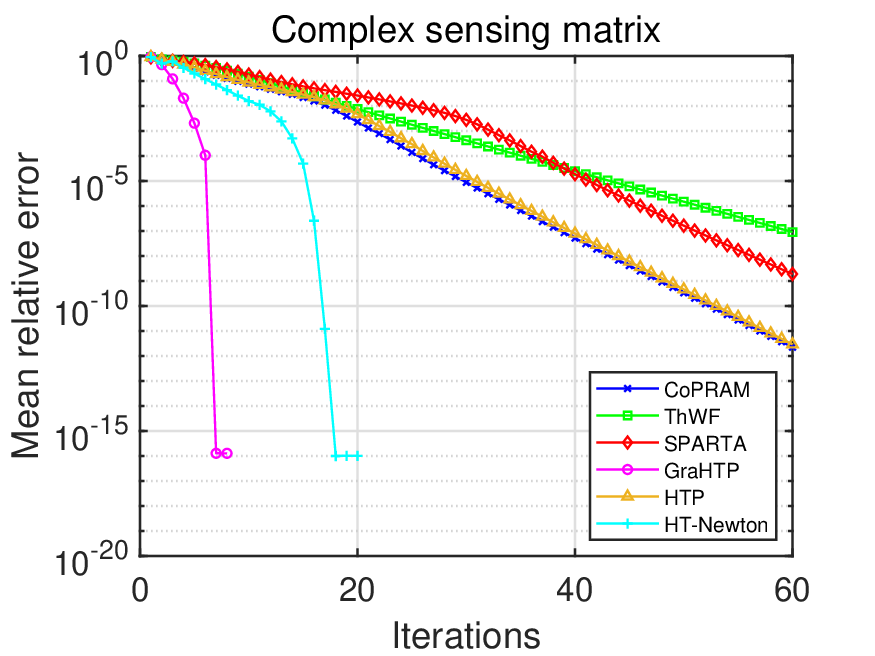}}\\
\subfigure[sparsity $s = 20$]{
\includegraphics[width=0.22\textwidth,trim=3pt 3pt 30pt 5pt,clip]{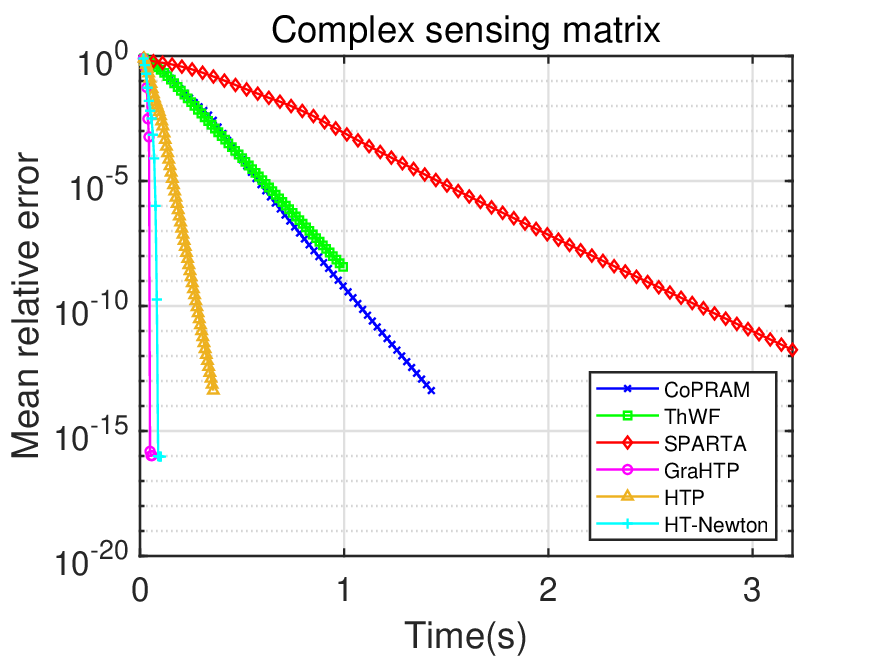}}
\subfigure[sparsity $s = 30$]{
\includegraphics[width=0.22\textwidth,trim=3pt 3pt 30pt 5pt,clip]{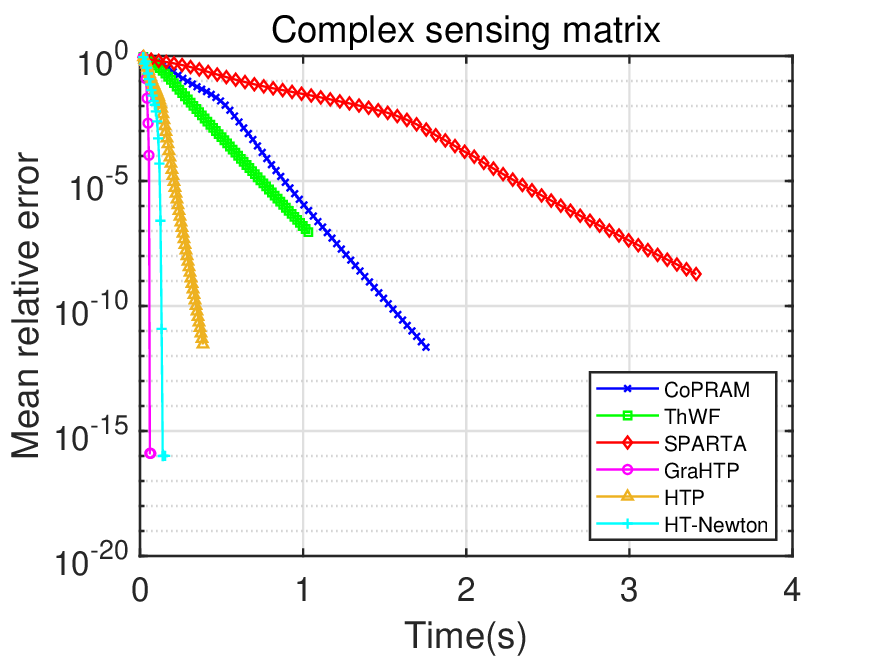}} 
\caption{Relative error versus number of iterations 
 ((a),(b)) and relative error versus running time ((c),(d)) for different algorithms, with $n = 3000$
and $m = 2000$. The results represent the average of $100$ independent trial runs.}
\label{figure1.1}
\end{center}
\end{figure}

\begin{figure}[ht!]
\begin{center}
\subfigure[sparsity $s = 20$]{
\includegraphics[width=0.22\textwidth,trim=3pt 3pt 30pt 5pt,clip]{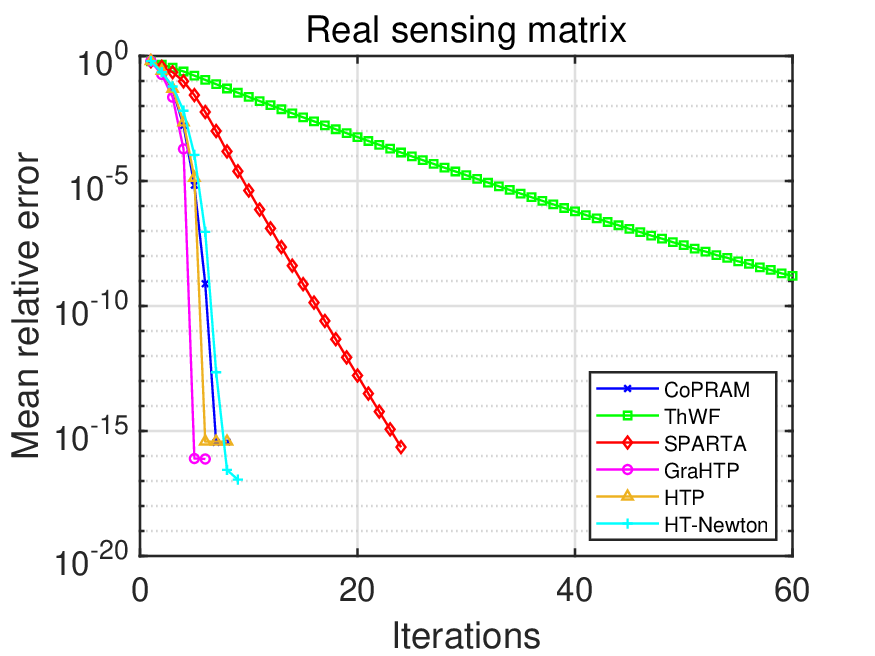}} 
\subfigure[sparsity $s = 30$]{
\includegraphics[width=0.22\textwidth,trim=3pt 3pt 30pt 5pt,clip]{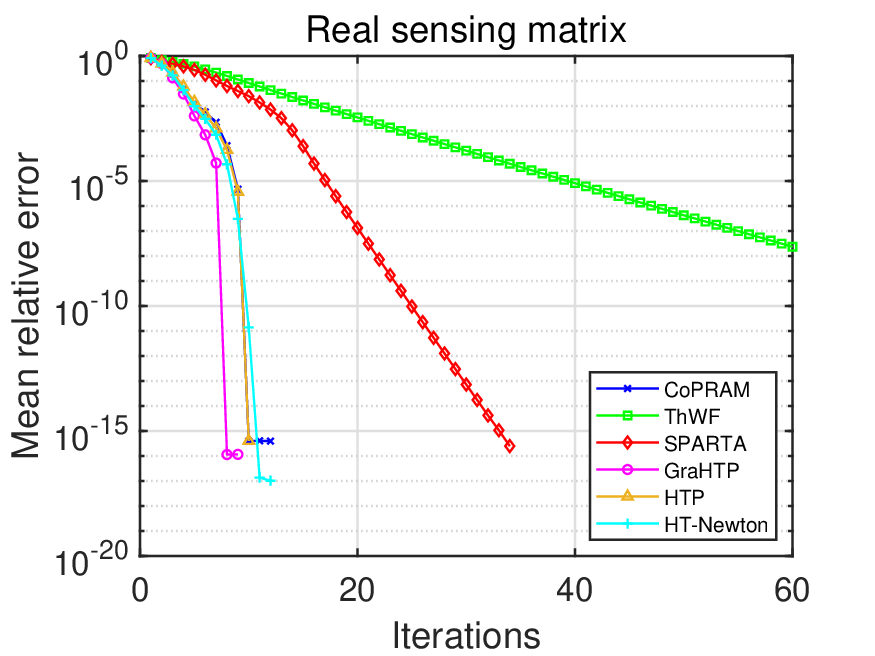}}\\
\subfigure[sparsity $s = 20$]{
\includegraphics[width=0.22\textwidth,trim=3pt 3pt 30pt 5pt,clip]{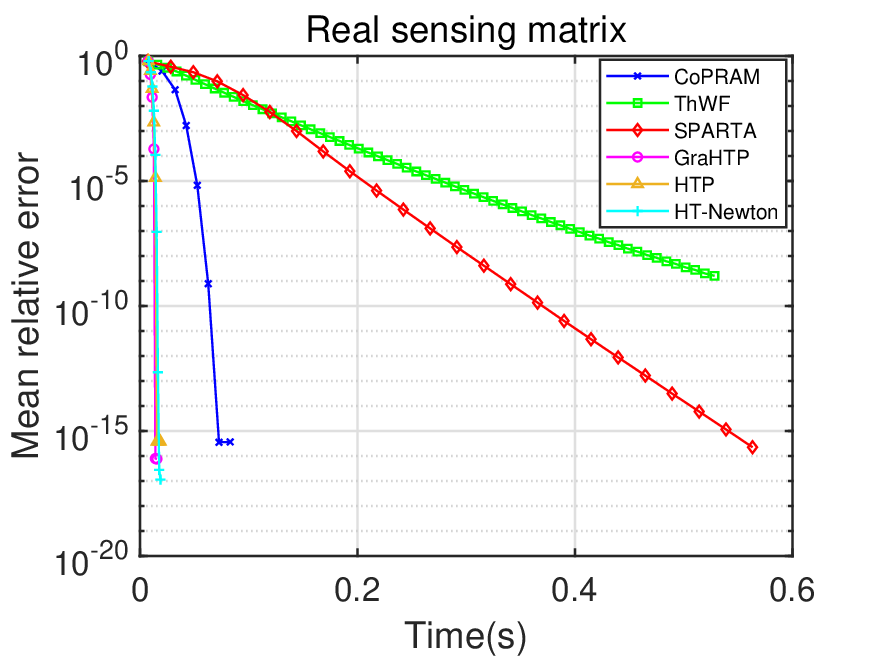}} 
\subfigure[sparsity $s = 30$]{
\includegraphics[width=0.22\textwidth,trim=3pt 3pt 30pt 5pt,clip]{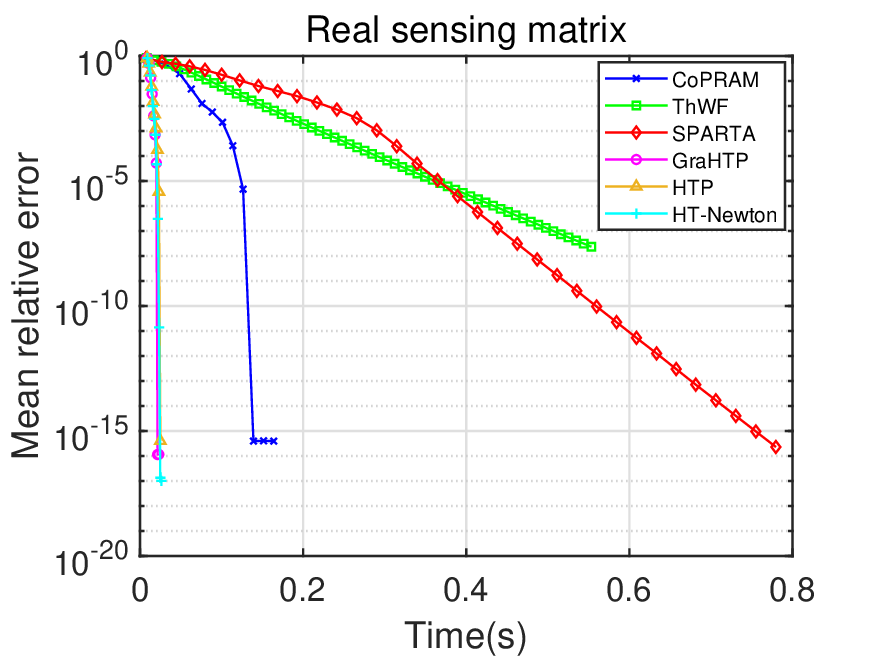}}
\caption{Relative error versus number of iterations 
 ((a),(b)) and relative error versus running time ((c),(d)) for different algorithms, with $n = 3000$
and $m = 2000$. The results represent the average of $100$ independent trial runs.}
\label{figure1.2}
\end{center}
\end{figure}

\textbf{Relative error.}
   We compared the convergence performance of different algorithms in the cases of complex and real sensing matrix, respectively, and the results are shown in \Cref{figure1.1} and \Cref{figure1.2}, namely. The maximum number of iterations is $60$. The signal length and the sample size are fixed at $n = 3000$ and $m = 2000$. The  sparsity levels are set to $s = 20, 30$.  

The results show that our algorithm achieves $r\left(\hat{\bm{x}}, \bm{x}^\dag\right) \leq 10^{-15}$ requiring fewer iterations and less time than all other methods when the sensing vectors are complex-valued. For real-valued sensing vectors, to achieve $r\left(\hat{\bm{x}}, \bm{x}^\dag\right) \leq 10^{-15}$, the computation time of our algorithm is comparable to that of HTP and HT-Newton.

\begin{figure}[ht!]
   \centering
\includegraphics[width=0.36\textwidth,trim=0pt 4pt 10pt 12pt,clip]{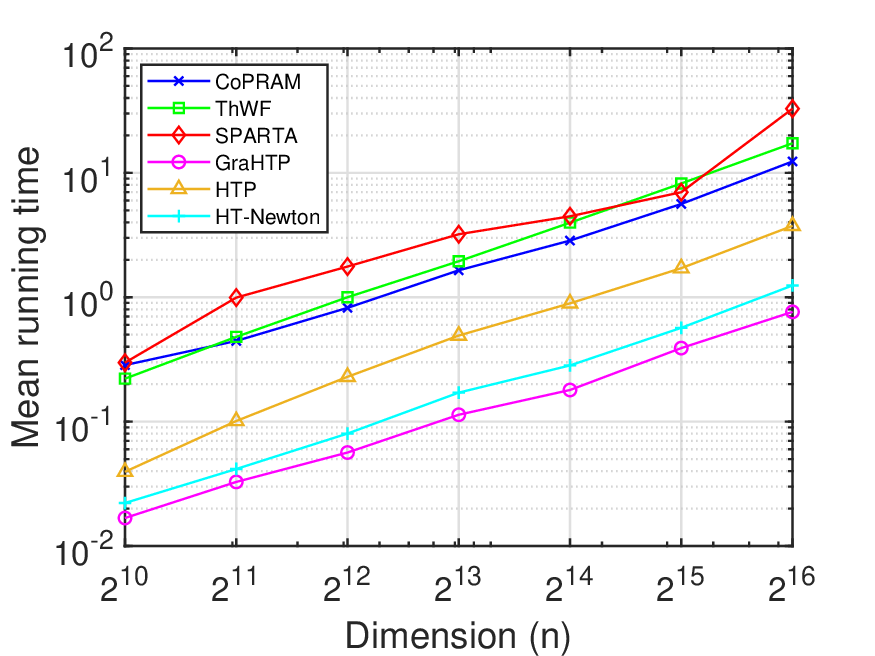}
    \caption{Running time for successful recovery versus signal dimension for different algorithms, with fixed sample size $m = 2120$ and sparsity $20$. All results were obtained by averaging 100 independent experiments with those fail trials filtered out.}
    \label{figure2}
\end{figure}

\textbf{Running time comparison.}
We compared several algorithms in terms of running time for successful recovery ($r\left(\hat{\bm{x}}, \bm{x}^\dag\right) \leq 10^{-6}$). In the experiment, the sample size $m$ is fixed at $2120$, the sparsity is fixed at $20$, and the signal length varies from $2^{10} \sim 2^{16}$. The maximum number of iterations is $60$. The results without those fail trials are shown on the \Cref{figure2}. According to the experiment, our algorithm GraHTP has a better efficiency performance than others.

 The computational complexities of CoPRAM, SPARTA and ThWF are both $\mathcal{O}(mn)$ per iteration. HTP and HT-Newton cost $\mathcal{O}(mn + s^2m)$ per iteration, and GraHTP costs $\mathcal{O}(mn + L_0s^2m)$ per iteration, which are of the comparable order as $\mathcal{O}(mn)$ in the case of $s \ll n$.

\begin{figure}[ht!]
   \centering
\includegraphics[width=0.23\textwidth,trim=7pt 2pt 20pt 5pt,clip]{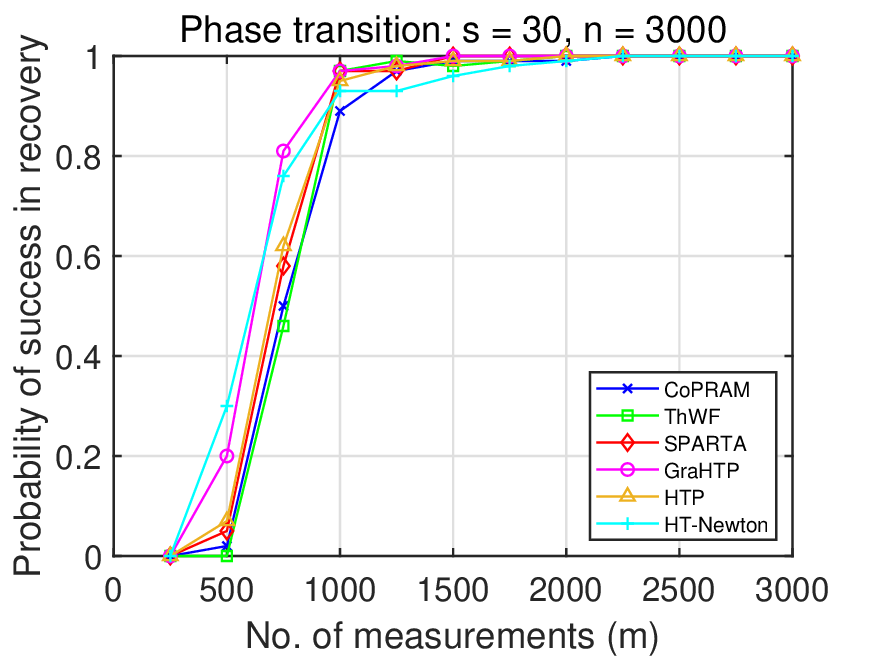}
\includegraphics[width=0.23\textwidth,trim=7pt 2pt 20pt 5pt,clip]{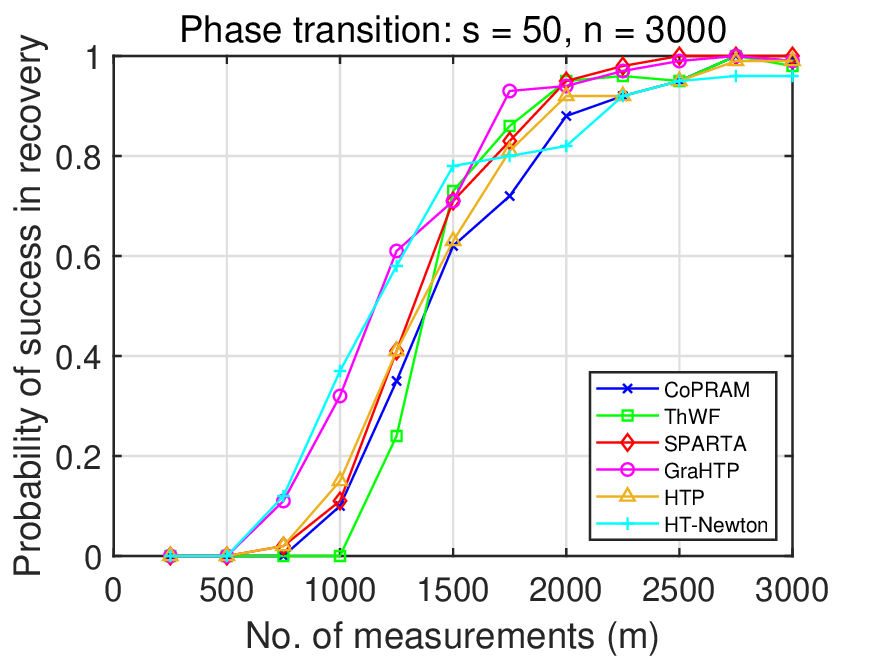}\\
    \caption{Phase transition for different algorithms, the results were obtained by averaging 100 independent experiments.}
    \label{figure3}
\end{figure}

\begin{figure}[ht!]
\centering   
\subfigure[ThWF]{\includegraphics[width=0.2\textwidth,trim=40pt 0pt 40pt 8pt,clip]{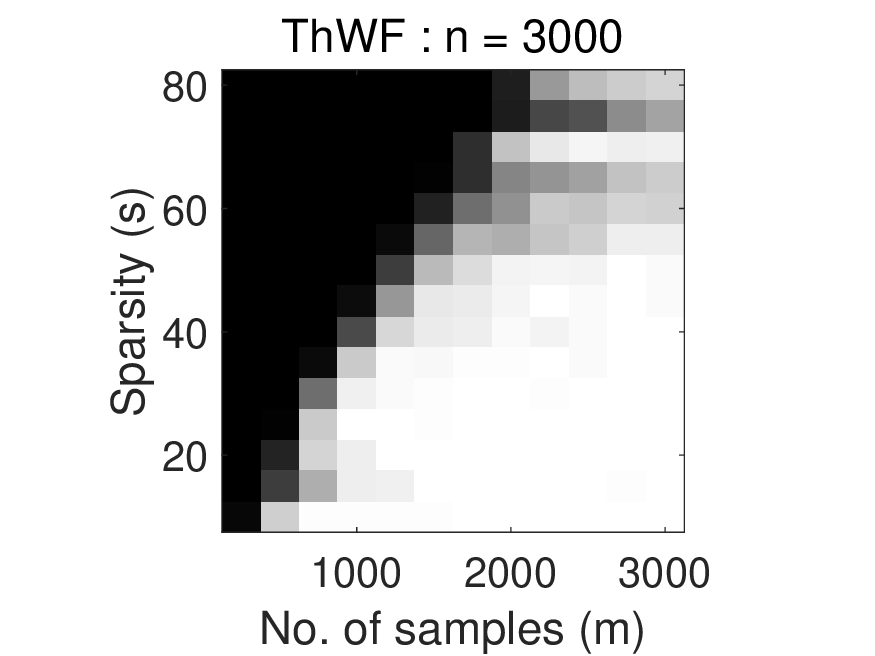}}
\subfigure[CoPRAM]{\includegraphics[width=0.2\textwidth,trim=40pt 0pt 40pt 8pt,clip]{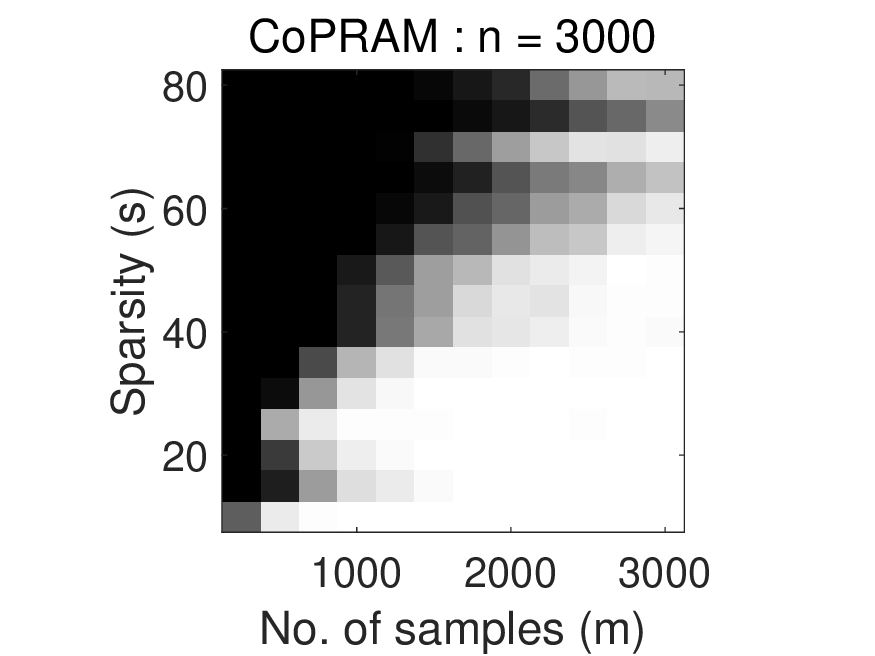}}\\
\subfigure[SPARTA]{\includegraphics[width=0.2\textwidth,trim=40pt 0pt 40pt 8pt,clip]{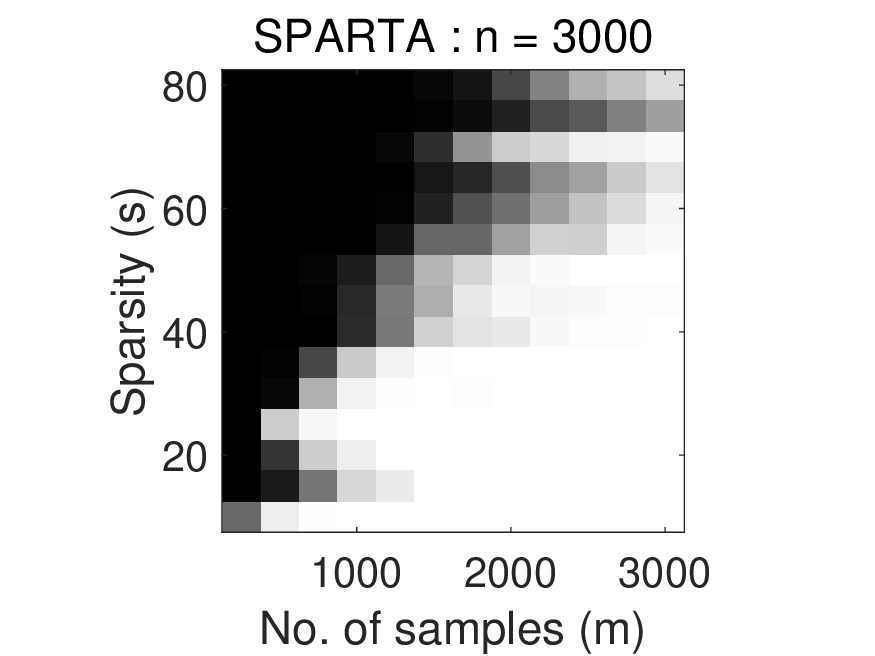}}
\subfigure[HTP]{\includegraphics[width=0.2\textwidth,trim=40pt 0pt 40pt 8pt,clip]{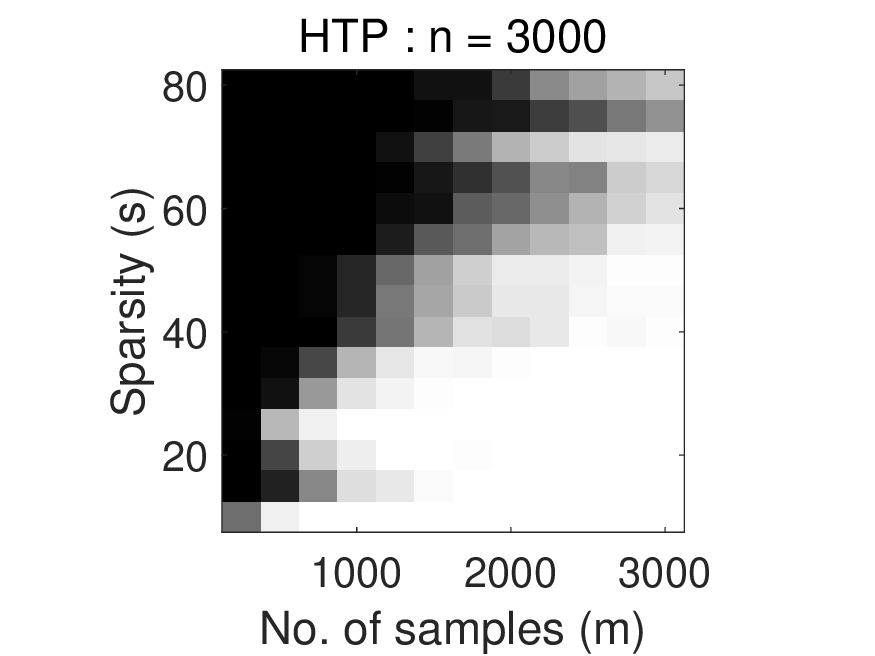}}\\
\subfigure[GraHTP]{\includegraphics[width=0.2\textwidth,trim=40pt 0pt 40pt 8pt,clip]{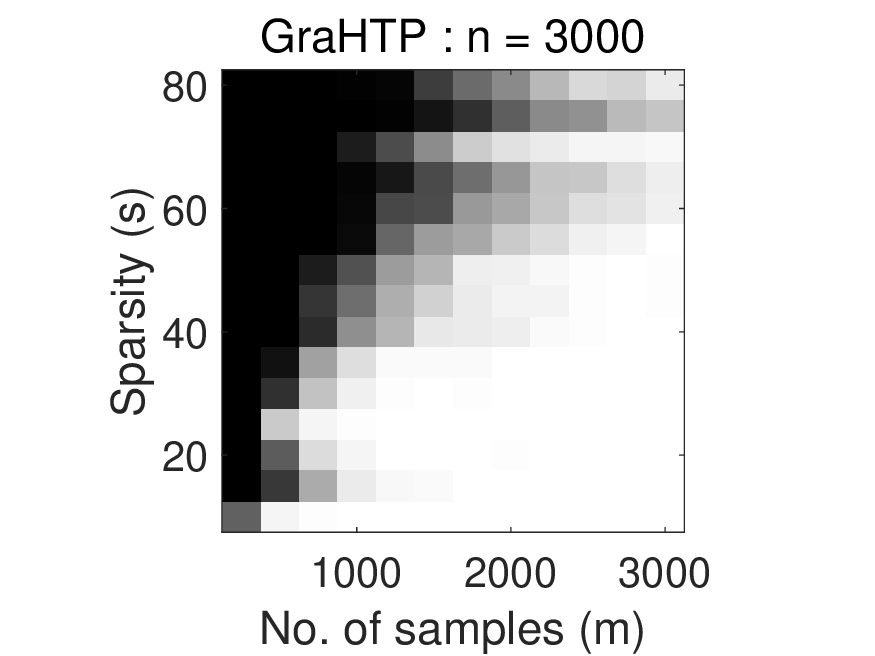}}
\subfigure[HT-Newton]{\includegraphics[width=0.2\textwidth,trim=40pt 0pt 40pt 8pt,clip]{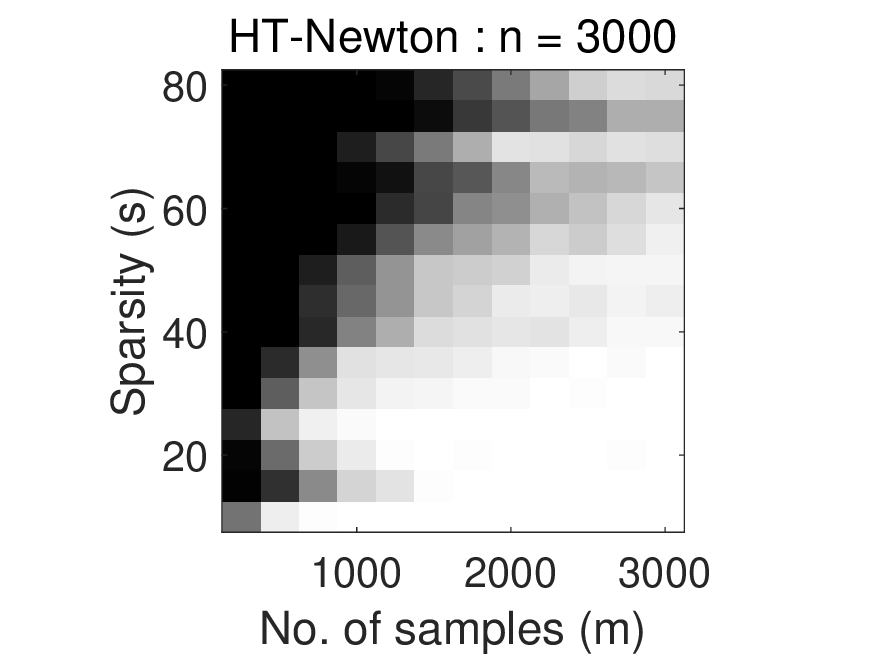}}
\caption{Phase transition for different algorithms, white block means $100\%$ successful recovery, black block means $0\%$ successful recovery and gray block means the rate of successful recovery between $0\%$ and $100\%$. The results were obtained by averaging 100 independent experiments.}
\label{figure4}
\end{figure}

\textbf{Phase transition.}
We present the results of comparing the recovery success rate ($r\left(\hat{\bm{x}}, \bm{x}^\dag\right) \leq 10^{-6}$) of different algorithms in \Cref{figure3,figure4}. The
maximum number of iterations is 100. In the first experiment, the signal length is fixed at $n = 3000$, the sparsity levels are set to $s = 30, 50$, and the sample size $m$ varies from $250\sim3000$. The result is shown in \Cref{figure3}. 
In the second experiment, the signal length is fixed at $n = 3000$, the sparsity $s$ varies from $10\sim80$ with grid size $5$, and the sample size $m$ varies from $250\sim3000$ with grid size $250$. The result is shown in \Cref{figure4}, the gray level of a block means the success recovery rate under the given sparsity and sample size. 

\begin{figure}[ht!]
   \centering
\includegraphics[width=0.3\textwidth,trim=2pt 2pt 20pt 5pt,clip]{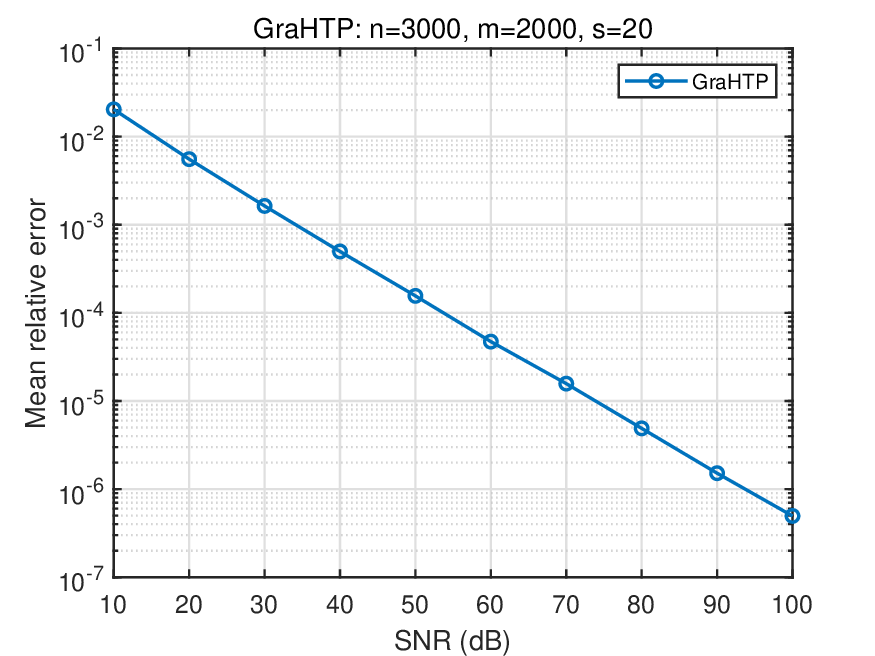}
    \caption{Relative error versus SNR for GraHTP, the results were obtained by averaging 100 independent experiments.}
    \label{figureRE}
\end{figure}

\textbf{Signal-to-Noise Ratio.}
We present the results of the convergence performance of our algorithm GraHTP under different cases of Signal-to-Noise Ratio (SNR) in \Cref{figureRE}. The SNR are set from $10$ to $100$ dB, the signal length is fixed to $3000$, the sample size is set as $2000$, and the sparsity is set as $20$. 

\begin{figure}[ht!]
 \centering
\subfigure[True Signal]{\includegraphics[width=0.22\textwidth,trim=0pt 50pt 0pt 0pt,clip]{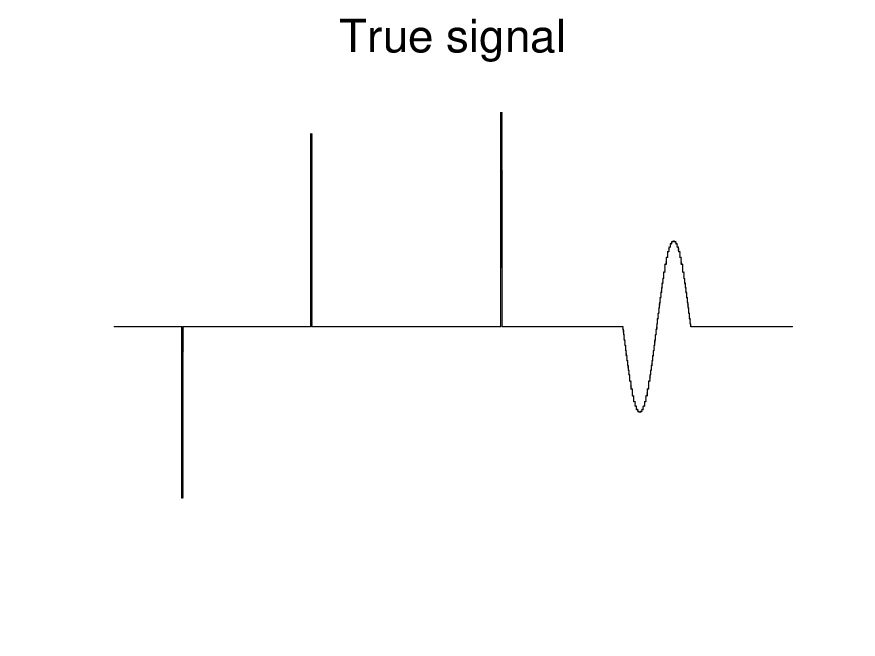}}\\
\subfigure[CoPRAM]{\includegraphics[width=0.22\textwidth,trim=0pt 50pt 0pt 0pt,clip]{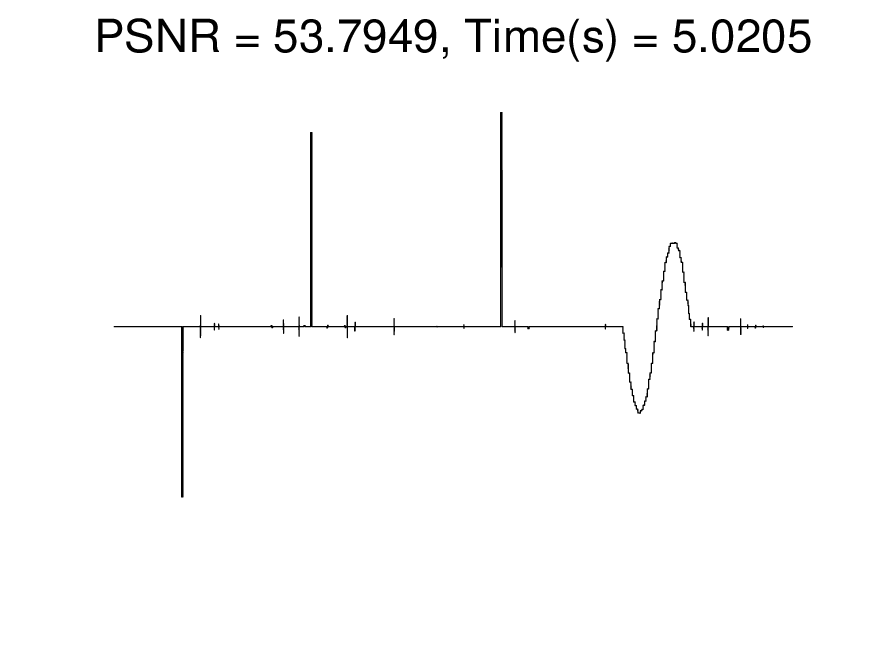}}
\subfigure[ThWF]
{\includegraphics[width=0.22\textwidth,trim=0pt 50pt 0pt 0pt,clip]{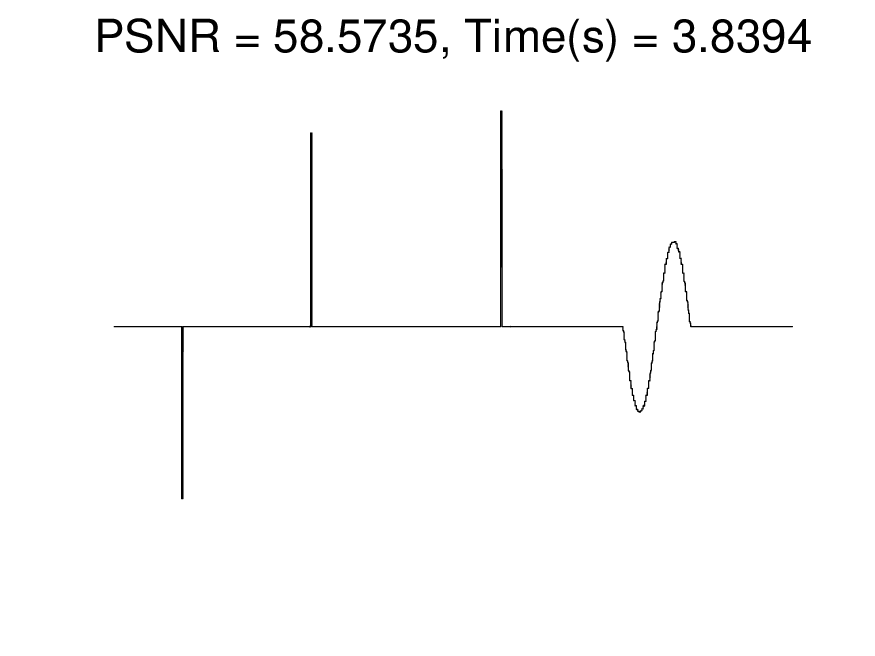}}\\
\subfigure[SPARTA]{\includegraphics[width=0.22\textwidth,trim=0pt 50pt 0pt 0pt,clip]{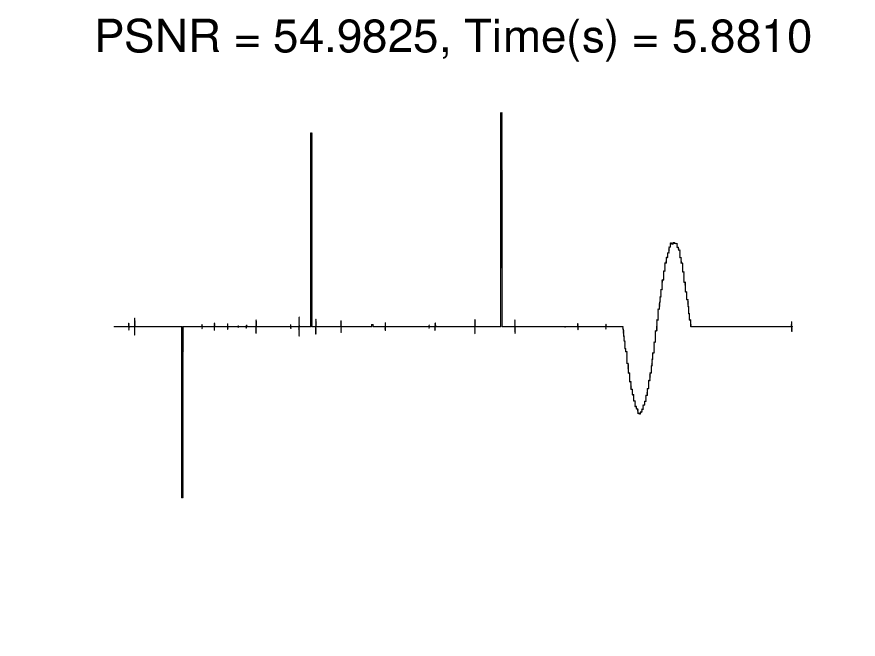}}
\subfigure[HTP]{\includegraphics[width=0.22\textwidth,trim=0pt 50pt 0pt 0pt,clip]{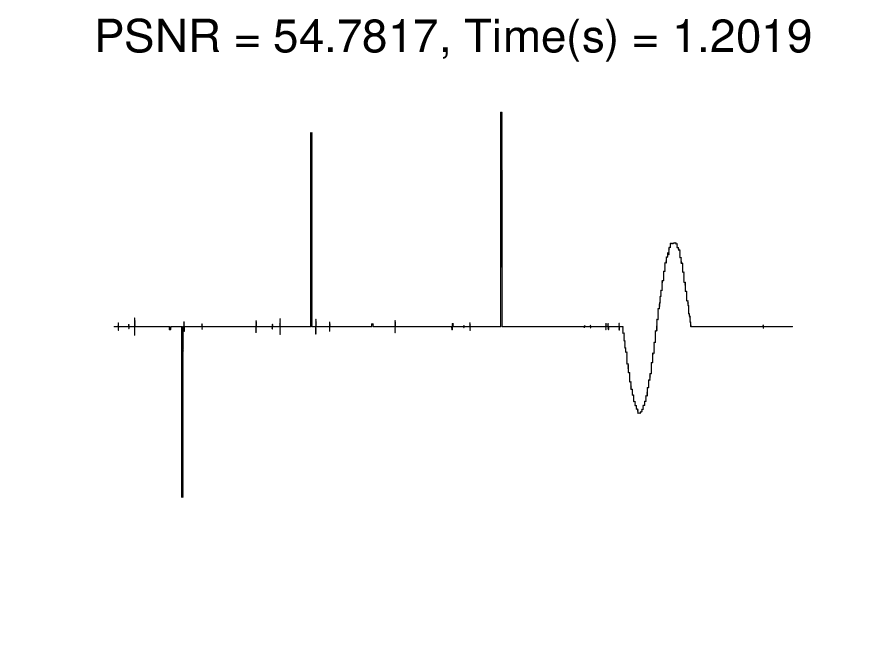}}\\
\subfigure[GraHTP]{\includegraphics[width=0.22\textwidth,trim=0pt 50pt 0pt 0pt,clip]{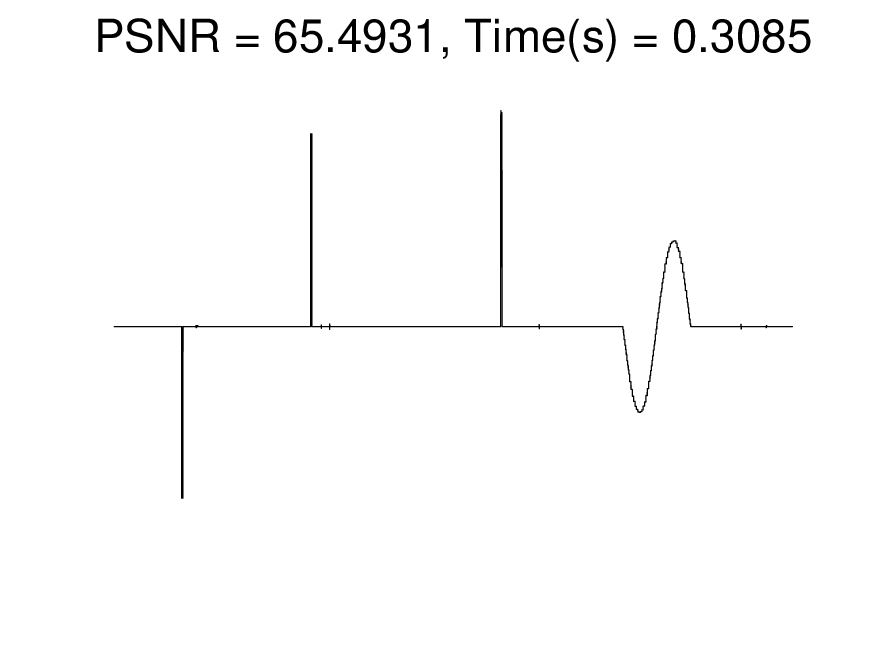}}
\subfigure[HT-Newton]{\includegraphics[width=0.22\textwidth,trim=0pt 50pt 0pt 0pt,clip]{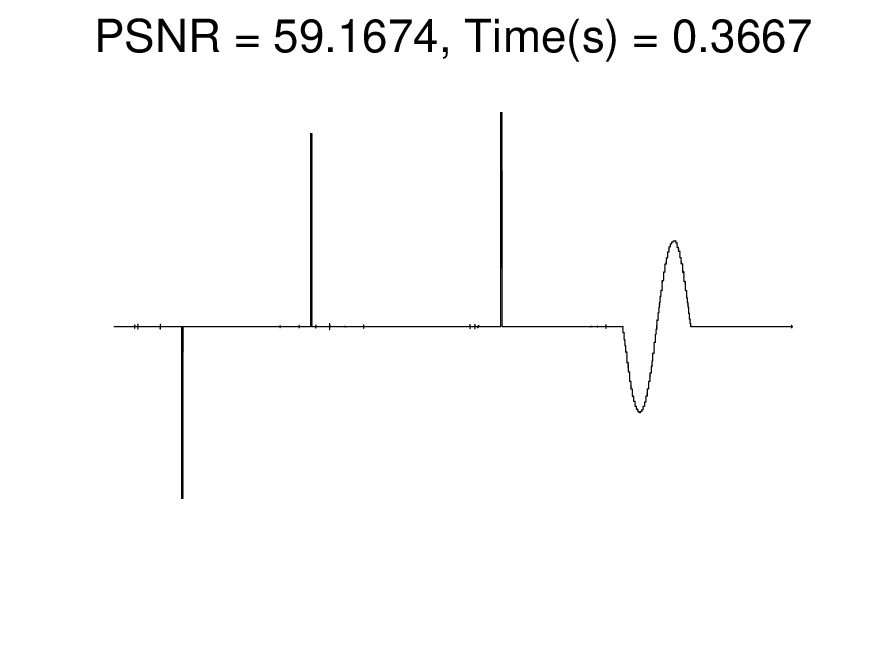}}
\caption{1-D signal reconstruction for different algorithms. PSNR represents peak signal-to-noise ratio and Time(s) is the running time in seconds.}
\label{figure5}
\end{figure}

\textbf{1-D signal reconstruction.}
Now we test the performance of different algorithms on recovering a 1-D signal from phaseless noisy measurements, the results of which are shown in \Cref{figure5}. The sampling matrix $\bm{A}$ is of size $2800 \times 8000$ and is constructed from
a complex random Gaussian matrix and an inverse wavelet transform (with four levels of Daubechies 1 wavelet). The noise level is $\sigma = 0.05$. The signal is sparse induced in the wavelet transformation and we set $s$ to $0.01n$ in the numerical experiment since the exact sparsity level is unknown in practical. The PSNR value is defined as
\begin{align*}
    \mathrm{PSNR}=10 \cdot \log \frac{\mathrm{V}^2}{\mathrm{MSE}}
\end{align*}
where $V$ represents the peak of the true signal and MSE is the mean squared error of the signal reconstruction. The result shows that our proposed algorithm GraHTP cost less time to achieve the higher PSNR in signal reconstruction.

\begin{figure}[ht!]
\centering 
\subfigure[True Image]{\includegraphics[width=0.23\textwidth,trim=47pt 50pt 47pt 16pt,clip]{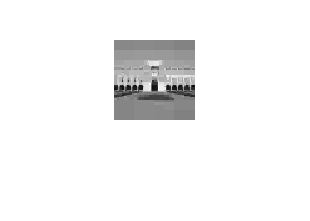}}\\
\subfigure[CoPRAM: PSNR=25.5369]{\includegraphics[width=0.23\textwidth,trim=47pt 50pt 47pt 16pt,clip]{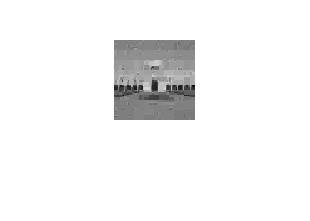}}
\subfigure[ThWF: PSNR=21.9714]{\includegraphics[width=0.23\textwidth,trim=47pt 50pt 47pt 16pt,clip]{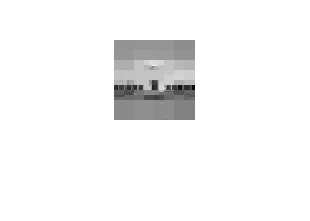}}\\
\subfigure[SPARTA: PSNR=30.8928]{\includegraphics[width=0.23\textwidth,trim=47pt 50pt 47pt 16pt,clip]{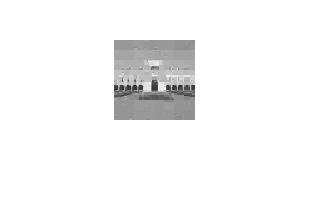}}
\subfigure[HTP: PSNR=30.9975]{\includegraphics[width=0.23\textwidth,trim=47pt 50pt 47pt 16pt,clip]{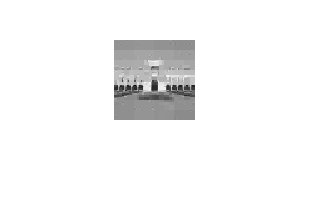}}\\
\subfigure[GraHTP: PSNR=44.2606]{\includegraphics[width=0.23\textwidth,trim=47pt 50pt 47pt 16pt,clip]{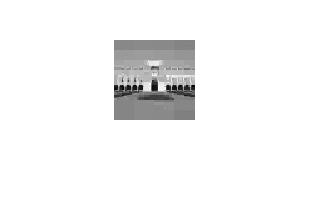}}
\subfigure[HT-Newton: PSNR=37.8018]
{\includegraphics[width=0.23\textwidth,trim=47pt 50pt 47pt 16pt,clip]{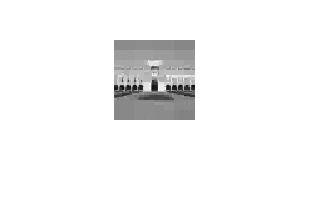}}
\caption{2-D image reconstruction for different algorithms, PSNR represents peak signal-to-noise ratio.}
\label{figure6}
\end{figure}

\textbf{2-D image reconstruction.}
\Cref{figure6} compares the performance of different algorithms in recovering a 2-D image (of size $64 \times 64$) with induced sparsity in the wavelet transform domain
from noisy phaseless measurements. We use a thresholded wavelet transform (with four level of Daubechies 1 wavelet) of this image as the target signal. The noise level in the measurements is set to $\sigma = 0.03$. The matrix $\bm{A}$ is a complex random Gaussian matrix of size $3600 \times 4096$. In the numerical experiment, the exact sparsity level is not exactly known, and the sparsity $s$ is set to $400$ (around $0.1n$) for image reconstruction. The results show that in this experiment the PSNR of our proposed algorithm GraHTP is higher than other methods.

\subsection{Complex-valued Signal Case}\label{complex signal case}
Next, we will give the results of the numerical experiments for complex-valued signal. 

\begin{figure}[ht!]
\centering
\subfigure[sparsity $s = 20$]{
\includegraphics[width=0.23\textwidth,trim=3pt 0pt 30pt 0pt,clip]{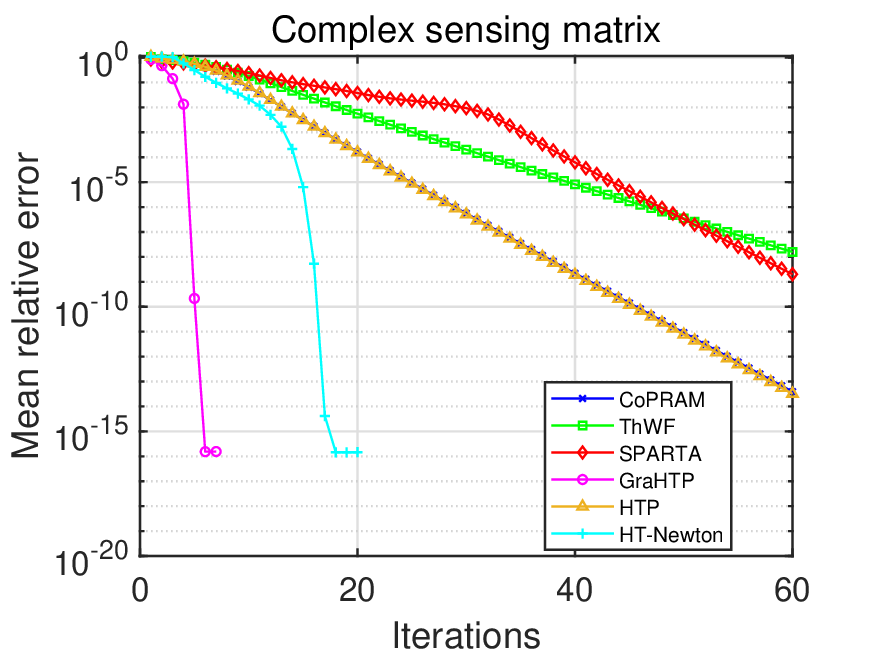}}
\subfigure[sparsity $s = 20$]{
\includegraphics[width=0.23\textwidth,trim=3pt 0pt 30pt 0pt,clip]{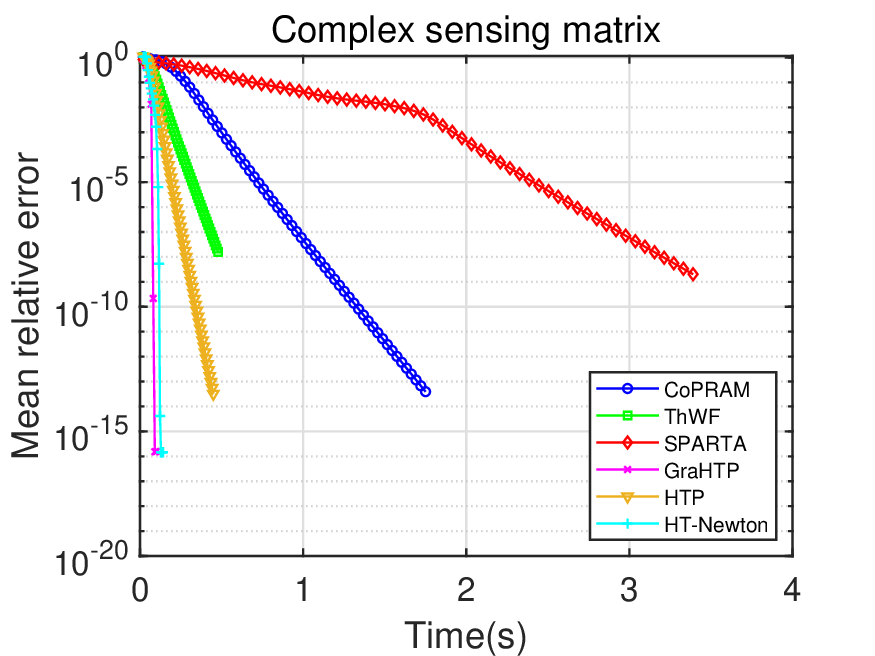}}
\caption{Relative error versus number of iterations (a) and relative error versus running time (b) for different algorithm, in the case of complex-valued signal with $n = 3000$ and $m = 2500$. The results represent the average of $100$ independent trial runs with those fail trials filtered out.}
\label{figure7}
\end{figure}

\textbf{Relative error.}
We compare different algorithms in the case where the true signal and sensing vectors are all complex-valued, the experimental results in \Cref{figure7} 
show the convergence of different algorithms. In this experiment, the signal length is fixed at $n = 3000$, the sample size is fixed at $m = 2500$ and the sparsity levels are set at $s = 20$. The maximum number of iterations is $60$. The results show that our algorithm GraHTP requires fewer iterations and less time while achieving $r(\hat{\bm{x}},\bm{x}^\dag) \leq 10^{-15}$.

\begin{figure}[ht!]
   \centering \includegraphics[width=0.3\textwidth,trim=3pt 0pt 10pt 0pt,clip]{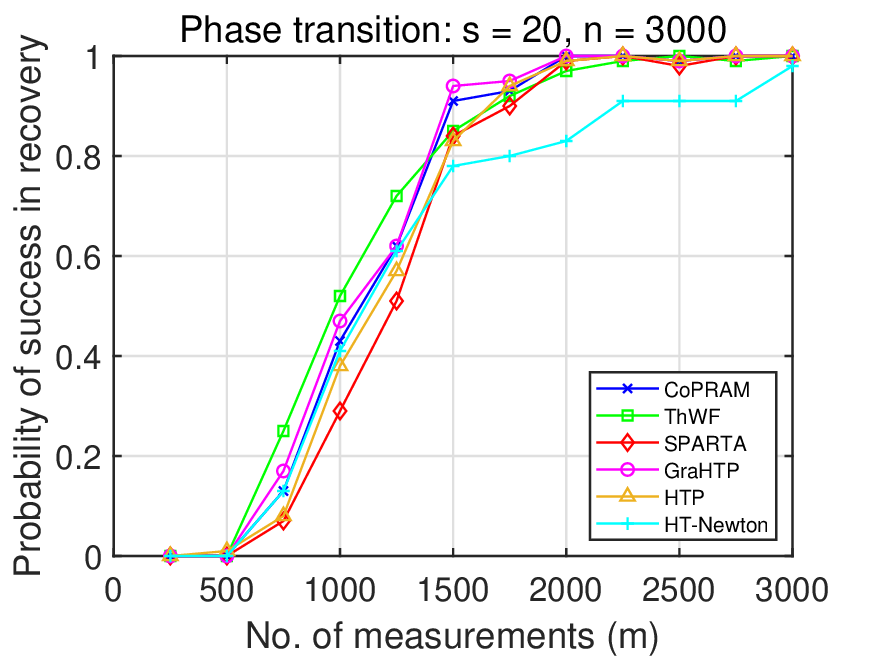}\\
    \caption{Phase transition for different algorithms. The results represent the average of $100$ independent trial runs.}
    \label{figure8}
\end{figure}

\textbf{Phase transition.}
We present the results of comparing the recovery success rate ($r(\hat{\bm{x}},\bm{x}^\dag) \leq 10^{-6}$) of several algorithms in \Cref{figure8}.  
In this experiment, the true signal dimension is fixed at $n = 3000$, the sparsity is fixed at $s = 20$, and the corresponding sample sizes $m$ vary from $250\sim3000$. The maximum number of iterations is 100. The $x$ axis in the figure represents the sample size $m$, and the $y$ axis represents the successful reconstruction rate. 

\subsection{Partial Discrete Fourier Transform Matrix}
\textbf{Relative error.}
   In this experiment, we compare the number of iterations required for our proposed algorithm GraHTP in the case where the sensing matrix $\bm{A}$ is a partial Discrete Fourier Transform (DFT) matrix, and the true signal is real-valued. The spectral initialization method used in this work is not directly applicable in the Fourier sensing setting since it does not have comparable theoretical guarantees as in the Gaussian sensing case. Therefore, to evaluate the local convergence behavior of GraHTP independently of initialization issues, we empirically used an initial estimate $\bm z^0$ satisfying the initialization condition $r(\bm z^0, \bm x^\dagger) \leq 0.8$ as the input of our algorithm. The signal dimension is fixed to be $n = 2000$, the sample size is fixed to be $m = 1500$, the sparsity of true signal is set to be $s = 20$ and $s = 30$, respectively and the maximum number of iterations is $10$. The $x$ axis in the figure represents the iterations, and the $y$ axis represents the relative error. We see that our proposed algorithm only need a few number of iterations for achieving $r\left(\hat{\bm{x}}, \bm{x}^\dag\right) \leq 10^{-15}$.
   
\begin{figure}[ht!]
\centering
\subfigure[sparsity $s = 20$]{
\includegraphics[width=0.23\textwidth,trim=5pt 0pt 30pt 5pt,clip]{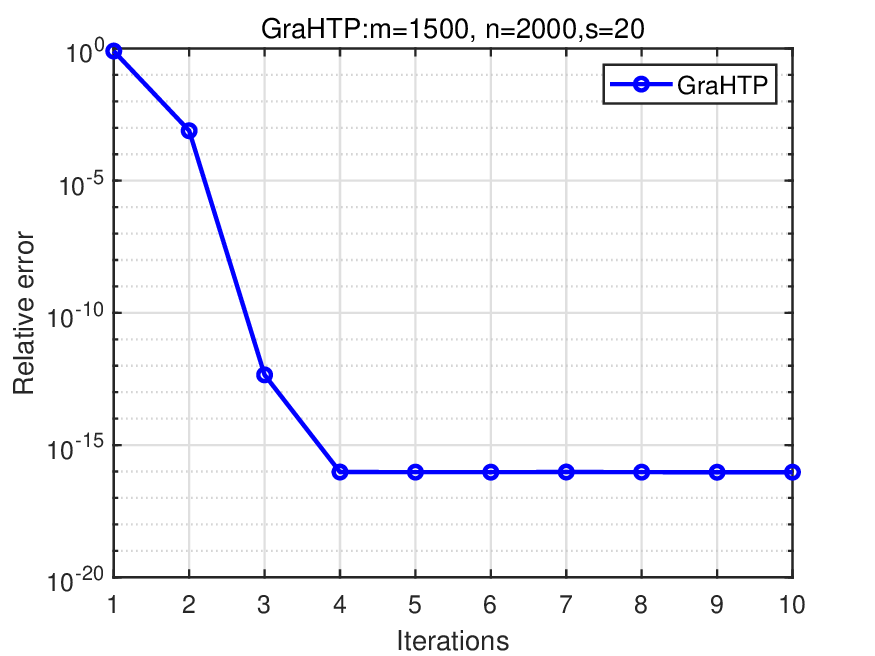}}
\subfigure[sparsity $s = 30$]{
\includegraphics[width=0.23\textwidth,trim=5pt 0pt 30pt 5pt,clip]{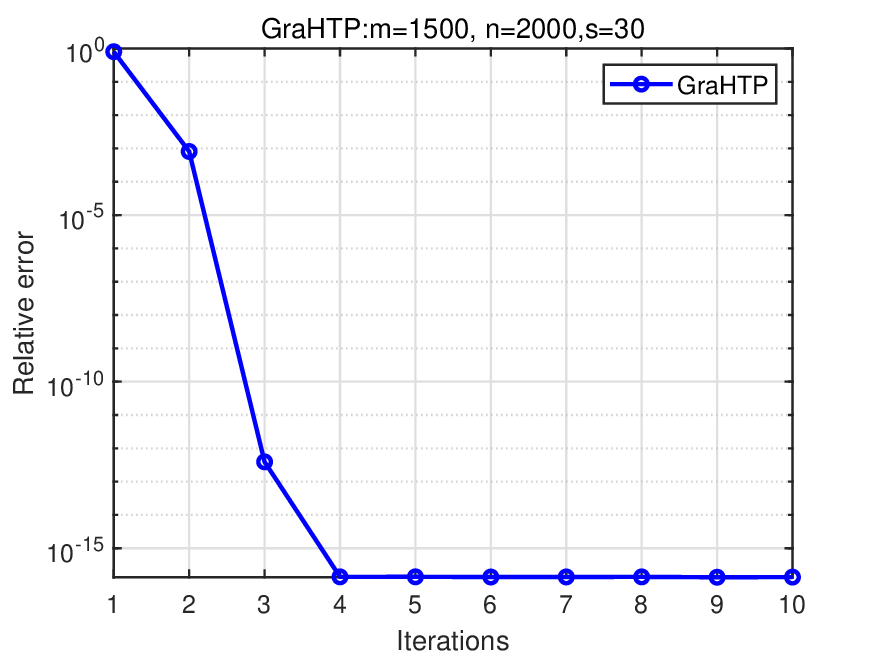}}
\caption{Partial DFT sensing matrix: Relative error versus number of iterations for GraHTP, with fixed signal dimension $n = 2000$, sample size $m = 1500$.}
\label{figure9}
\end{figure}

\noindent{\bf Acknowledgements.} This work is supported by the National Key Research and Development Program of China (Grant No. 2025YFA1017800), the National Science Foundation of China (No. 12371424, No. 12401121, No. U24A2002, No. 12561160122), the Natural Science Foundation of Hubei Province (No. 2024AFE006), the Science and Technology Innovation Platform Program of Hubei Province (No. 2025CSA119) , the Fundamental Research Funds for the Central Universities(No. 2042026fg0009), the Shenzhen Municipal Stability Support Program Project (No. 8960117/0123).

\section{Appendix}\label{proofs}
\subsection{Key Lemmas}\label{proof1}
In this subsection, we first give several key lemmas, which are essential for proving theoretical results \Cref{local convergence}, \Cref{noisy case} and \Cref{global}.

\begin{lemma}
	\label{S-ES}
Let $\bm{b}_{j}, \bm{c}_{j} \in \R^n , j=1,\cdots,m$ be i.i.d. $\mathcal{ { N }}(\mathbf{0},\bm{I}_n/2)$ random vectors. For any sparse vector $\bm{z} \in \R^n$ independent from $\{\bm{b}_{j}, \bm{c}_{j}\}_{j=1}^m$ and for any subset $\T \subseteq [n]$ with $|\T | \leq t$ ($t$ is an integer less than $n$) such that  $\supp(\bm{z}) \subseteq \T$, denote

\begin{align*}
\bm{\Phi}~&:=~\frac{1}{m} \sum_{j=1}^m\big|{\bm{b}}_{j,\T}^{\top} \bm{z}_{\T}\big|^2 {\bm{b}}_{j,\T} {\bm{b}}_{j,\T}^{\top},
\end{align*}
and
\begin{align*}
\bm{\Psi}~&:=~\frac{1}{m} \sum_{j=1}^m\big({\bm{b}}_{j,\T}^{\top} \bm{z}_{\T}\big)\big( {\bm{c}}_{j,\T}^{\top} \bm{z}_{\T}\big){\bm{b}}_{j,\T} {\bm{c}}_{j,\T}^{\top} .
\end{align*}
Then, for any $\delta  \in (0,1)$,  with probability at least  $1-c_1 \delta^{-2} m^{-1} - m^{-4}-c_2 \exp \left(-c_3 \delta^2 m / \log m\right)$,  it holds that
\begin{equation}
\begin{aligned}\label{Phi-EPhi}
\norm{\bm{\Phi}- \frac{1}{4}(\|\bm{z}_{\T} \|_2^2\bm{I}_{|\T|} + 2\bm{z}_{\T} \bm{z}_{\T}^{\top} )} ~\leq~ 
\frac{\delta}{4}\norm{\bm{z} }^2, \cr
\norm{\bm{\Psi}- \frac{1}{4}\bm{z}_{\T} \bm{z}_{\T}^{\top} } ~\leq~ \frac{\delta}{4}\norm{\bm{z} }^2,
\end{aligned}   
\end{equation}
provided $m \geq C(\delta) t \log (n/t)$, where $C(\delta)$ is a constant depending on $\delta$, and $c_1$, $c_2$ and $c_3$ are positive absolute constants.
\end{lemma}

\begin{IEEEproof}
The proof of this lemma is similar to that of Lemma B.1 in \cite{cai2024a} and Lemma 21 in \cite{sun2018geometric}.
\end{IEEEproof}

The following lemma is an extension of Lemma V.5 in \cite{gao2017phaseless}. 
For the ease of notations, we denote
\begin{align}\label{def:H}
\bm{H}(\bm{z})
& : = ~\bm{J}(\bm{z})^{\top}\bm{J}(\bm{z})  \cr
&= ~\frac{1}{m} \sum_{j=1}^m\Big((\bm{a}_{jR }^{\top} \bm{z})^2 \bm{a}_{jR } \bm{a}_{j R}^{\top}+(\bm{a}_{jI }^{\top} \bm{z})^2 \bm{a}_{jI } \bm{a}_{jI }^{\top} \cr
& \quad +~(\bm{a}_{jR }^{\top} \bm{z})(\bm{a}_{jI }^{\top} \bm{z})(\bm{a}_{jI } \bm{a}_{jR }^{\top}+\bm{a}_{jR } \bm{a}_{jI }^{\top})\Big).
\end{align}

\begin{lemma}\label{PHI}
Let $\bm{a}_{j} \in \C^n , j=1,\cdots,m$ be i.i.d. $\mathcal{ {C N }}(\mathbf{0},\bm{I}_n)$ Gaussian random vectors. For any sparse vector $\bm{z} \in \R^n$  and subset $\T \subseteq [n]$ such that $\supp(\bm{z}) \subseteq \T$ with $|\T | \leq t$ ($t$ is an integer less than $n$), under the event \eqref{Phi-EPhi}, 
we have
\begin{align}\label{J'J}
(\frac{1}{2}-\delta)\norm{\bm{z}}^2 ~\leq~ 
\norm{ \bm{H}_{\T,\T}(\bm{z})} ~\leq~ (2+\delta)\norm{\bm{z}}^2,
\end{align} 
and $\bm{H}_{\T,\T}(\bm{z})$ is invertible and
\begin{align}
\norm{\big(\bm{H}_{\T,\T}(\bm{z})\big)^{-1}} ~\leq~ \frac{2}{(1-2\delta)\norm{\bm{z}}^2}.
\end{align}
\end{lemma}

\begin{IEEEproof}
Note that $\supp(\bm{z}) \subseteq \T$, which implies that $\bm{a}_{jR}^{\top}\bm{z} = \bm{a}_{jR,\T}^{\top}\bm{z}_{\T}$ and $\bm{a}_{jI}^{\top}\bm{z} = \bm{a}_{jI,\T}^{\top}\bm{z}_{\T}$, then according to \eqref{def:H} we have
\begin{align*}
&\bm{H}_{\T,\T}(\bm{z}) \cr
 =~&\frac{1}{m} \sum_{j=1}^m\Big((\bm{a}_{jR,\T }^{\top} \bm{z}_{\T})^2 \bm{a}_{jR,\T } \bm{a}_{jR,\T }^{\top}+(\bm{a}_{jI,\T }^{\top} \bm{z}_{\T})^2 \bm{a}_{jI,\T }\bm{a}_{jI,\T }^{\top} \cr
&+(\bm{a}_{jR,\T }^{\top} \bm{z}_{\T})(\bm{a}_{jI,\T }^{\top} \bm{z}_{\T})(\bm{a}_{jI,\T }\bm{a}_{jR,\T }^{\top}+\bm{a}_{jR,\T } \bm{a}_{jI,\T }^{\top})\Big) .
\end{align*}
Assuming event \eqref{Phi-EPhi}. 
By direct calculation, we have
\begin{align*}
\E\big(\bm{H}_{\T,\T}(\bm{z})\big) ~=~ \frac{1}{2} \big(\norm{\bm{z}}^2 \bm{I}_{|\T|} + 3\bm{z}_{\T} \bm{z}_{\T}^{\top}\big),
\end{align*}
and the minimum and maximum eigenvalues of $\E(\bm{H}_{\T,\T}(\bm{z}))$ 
\begin{align}
\label{lambda min}	
\lambda_{\min }\big(\E\big(\bm{H}_{\T,\T}(\bm{z})\big)\big) ~&=~\frac{1}{2}\norm{\bm{z}}^2, \\
\label{lambda max}	
\lambda_{\max }\big(\E\big(\bm{H}_{\T,\T}(\bm{z})\big)\big) ~&=~ 2\norm{\bm{z}}^2 .
\end{align}	
Then by \Cref{S-ES} we have 
\begin{align*}
&\norm{\bm{H}_{\T,\T}(\bm{z}) - \E\big(\bm{H}_{\T,\T}(\bm{z})\big) } \cr 
\leq ~&2\bigg\|\frac{1}{m} \sum_{j=1}^m(\bm{a}_{jR,\T }^{\top} \bm{z}_{\T})^2 \bm{a}_{jR,\T } \bm{a}_{jR,\T }^{\top} -  \frac{1}{4} (\norm{\bm{z}}^2 \bm{I}_{|\T|} + 2\bm{z}_{\T} \bm{z}_{\T}^{\top})\bigg\|_2 \cr
 &+  2\bigg\| \frac{1}{m} \sum_{j=1}^m(\bm{a}_{jR,\T }^{\top} \bm{z}_{\T})(\bm{a}_{jI,\T }^{\top} \bm{z}_{\T})\bm{a}_{jI,\T } \bm{a}_{jR,\T }^{\top} - \frac{1}{4}\bm{z}_{\T} \bm{z}_{\T}^{\top}\bigg\|_2 \cr
\leq ~&\delta \norm{\bm{z}}^2 ,
\end{align*}
thus, by Weyl's Perturbation Theorem (Corollary III.2.6 in \cite{bhatia1997matrix}) we obtain
\begin{align*}
\big|\lambda_{\max }\big(\bm{H}_{\T,\T}(\bm{z})\big) - \lambda_{\max }\big(\E\big(\bm{H}_{\T,\T}(\bm{z})\big)\big)\big| ~&\leq~ \delta\norm{\bm{z}}^2,\cr
\big|\lambda_{\min }\big(\bm{H}_{\T,\T}(\bm{z})\big) -\lambda_{\min }\big(\E\big(\bm{H}_{\T,\T}(\bm{z})\big)\big)\big| ~&\leq~ \delta\norm{\bm{z}}^2.
\end{align*}
Together with \eqref{lambda min} and \eqref{lambda max}, we have
 \begin{align*}
(\frac{1}{2}-\delta)\norm{\bm{z}}^2 
~&\leq~	\lambda_{\min }\big(\bm{H}_{\T,\T}(\bm{z})\big)  ~\leq~ (\frac{1}{2}+\delta)\norm{\bm{z}}^2 
\cr (2 - \delta)\norm{\bm{z}}^2 ~&\leq~ \lambda_{\max }\big(\bm{H}_{\T,\T}(\bm{z})\big) ~\leq~ (2 + \delta)\norm{\bm{z}}^2 ,
\end{align*}
and
\begin{align*}
\lambda_{\max }\big(\big(\bm{H}_{\T,\T}(\bm{z})^{-1}\big)\big) ~=~ 1 / \lambda_{\min }\big(\bm{H}_{\T,\T}(\bm{z})\big) ~\leq~ \frac{2}{(1-2\delta)\norm{\bm{z}}^2},
\end{align*}
which completes the proof.
\end{IEEEproof}

\begin{lemma}\label{S-Sc}
Let $\bm{a}_{j} \in \C^n , j=1,\cdots,m$ be i.i.d. $\mathcal{ {C N }}(\mathbf{0},\bm{I}_n)$ Gaussian random vectors. 
 For any sparse vector $\bm{z} \in \R^n$ and  subsets $\cS, \T \subseteq [n]$ satisfying $\cS \subseteq \T$ with $|\T | \leq t$ ($t$ is an integer less than $n$), under the event \eqref{Phi-EPhi} 
 with $\supp(\bm{z}) \subseteq \cS$, it holds that
\begin{align*}
\norm{\bm{H}_{\cS,{\T\backslash\cS}}(\bm{z})} ~\leq~ \delta\norm{\bm{z}}^2 .
\end{align*}
\end{lemma}

\begin{IEEEproof}
Note that $\cS$ and $\T\backslash\cS$ are two disjoint subsets of $\T$, and we have
\begin{align*}
 \mathbb{E}[\bm{H}_{\cS,{\T\backslash\cS}}(\bm{z})]~=~ 
\bm{z}_{\cS}\bm{z}_{\T\backslash\cS}^{\top} ~=~ \bm{0},
\end{align*}
thus $\bm{H}_{\cS,{\T\backslash\cS}}(\bm{z})$ is a sub-matrix of $\bm{H}_{\T,\T}(\bm{z})-\E\big(\bm{H}_{\T,\T}(\bm{z})\big)$. According to \Cref{PHI} we then have
\begin{align*}	
\norm{\bm{H}_{\cS,{\T\backslash\cS}}(\bm{z})} 
~\leq~ \norm{\bm{H}_{\T,\T}(\bm{z})-\E\big(\bm{H}_{\T,\T}(\bm{z})\big)} ~\leq~ \delta\norm{\bm{z}}^2,
\end{align*}
which completes the proof.
\end{IEEEproof}

The following Lemma provides an upper bound for hard thresholding.
\begin{lemma}
\label{uv}
For any sparse vector $\bm{x} \in \R^n$ satisfying $\|\bm{x} \|_0 \leq s$ and any vector  $\bm{v} \in \R^n$, define $\bm{u} := \mathcal{H}_{s}(\bm{v})$, $\cS_{u} := \supp(\bm{u})$ and $\cS := \supp(\bm{x})$, then we have the following inequality: 
\begin{align*}
\|\bm{u} - \bm{x}\|_2^2 ~\leq~ \frac{3+\sqrt{5}}{2} \big\|\bm{v}_{{\cS_u}\cup \cS} - \bm{x}_{{\cS_u} \cup \cS}\big\|_2^2 .
\end{align*}
\end{lemma}

\begin{IEEEproof}
The proof of this lemma is a direct application of Theorem 1 in  \cite{shen2018tight} since $n \geq 2s$.
\end{IEEEproof}

The next lemma is about the upper bound on the estimation error of the vector obtained by applying the hard thresholding operator to the output of gradient descent.
 
\begin{lemma}\label{uk}
Given a sparse vector $\bm z^k$ that satisfies $\|\bm z^k -\bm x^\dag \|_2 \leq \gamma \|\bm x^\dag \|_2$ and $\|\bm z^k\|_0 \leq s$. Let $ \bm{a}_j\in \C^n $, $ j=1,\ldots, m $ be i.i.d. $\mathcal{CN}(0,\bm{I}_n)$ Gaussian random vectors independent with $\bm z^k$ and $\bm x^\dag$. The vector obtained by the projected gradient descent at $\bm z^k$ is defined as
\[
\bm u^k := \mathcal{H}_s(\bm{z}^k-\mu^k \nabla f(\bm{z}^k))
\]
with step size $\mu^k$, then under the event \eqref{Phi-EPhi}, 
it holds 
\begin{align}
\big\|\bm u^k - \bm x^\dag \big\|_2 
~\leq~ \rho\big\|\bm{z}^k-\bm{x}^\dag\big\|_2 ,
\end{align}
where $\rho = \sqrt{\frac{3+\sqrt{5}}{2}}\big( 1-\mu^{k}\norm{\bm{x}^\dag}^2(1 -2\delta  - (4 - 3\delta)\gamma - (1 + 3\delta)\gamma^2)\big)$ and $\mu^k \in \big(\frac{1}{(4-2\delta)(1-\gamma)^2\norm{\bm{x}^\dag}^2}, \frac{2}{(5+4\delta)(1+\gamma)^2\norm{\bm{x}^\dag}^2}\big)$. 
\end{lemma}

\begin{IEEEproof}\label{u_{k+1}}
Denote $\cS^{\dag} := \supp(\bm{x}^\dag)$, $\cS_{k} := \supp(\bm{z}^k)$, $\cS_{k+1} := \supp(\bm{u}^{k})$, $\T_{k+1} :=  \cS_{k}\bigcup\cS_{k+1} \bigcup \cS^{\dag}$, $\bm{h}^{k} := \bm{z}^{k} - \bm{x}^{\dag}$ and 
\[
\bm{v}^{k} ~:=~ \bm{z}^k-\mu^k \nabla f(\bm{z}^k),
\]
Then $\bm{u}^{k} =  \mathcal{H}_s(\bm{v}^{k})$.
According to the definition of $\T_{k+1}$,  the size of  $\T_{k+1}$ is at most $3s$, \textit{i.e.}, $|\T_{k+1}| \leq 3s$. 
Noting that $\bm{u}^{k}$ is the best $s$-term approximation of $\bm{v}^k$, then by \Cref{uv} and $\cS_{k+1} \bigcup \cS^{\dag} \subseteq \T_{k+1}$, we obtain
\begin{align}\label{u-x-tk+1}
\big\|\bm{u}^{k}-\bm{x}^\dag\big\|_2 
~&=~ \big\|\bm{u}^{k}_{\cS_{k+1} \cup \cS^{\dag}} - \bm{x}^\dag_{\cS_{k+1} \cup \cS^{\dag}} \big\|_2 \cr
~&\leq~  \sqrt{\frac{3+\sqrt{5}}{2}} \big\|\bm{v}^{k}_ {\T_{k+1}}-\bm{x}^\dag_{\T_{k+1}}\big\|_2.
\end{align}
Since $\cS_k, \cS^{\dag} \subseteq \T_{k+1}$, it holds that 
\begin{align*}
&~\nabla {f}_{\T_{k+1}}(\bm{z}^k)  \cr
= &~\frac{1}{m}\sum_{j=1}^m \Big((\bm{a}_{jR}^{\top}\bm{z}^k)^2 + (\bm{a}_{jI}^{\top}\bm{z}^k)^2 - (\bm{a}_{jR}^{\top}\bm{x}^\dag)^2 - 
(\bm{a}_{jI}^{\top}\bm{x}^\dag)^2 \Big)   \cr
&~ \cdot(\bm{a}_{jR,\T_{k+1}}\bm{a}_{jR}^{\top}\bm{z}^k + \bm{a}_{jI,\T_{k+1}}\bm{a}_{jI}^{\top}\bm{z}^k)  \cr
= &~\bm{J}_{\T_{k+1}}(\bm{z}^k)^{\top}\Big(\bm{J}_{\T_{k+1}}(\bm{z}^k ) + \bm{J}_{\T_{k+1}}(\bm{x}^\dag)\Big)(\bm{z}^k_{\T_{k+1}} - \bm{x}^\dag_{\T_{k+1}}) \cr
= &~ 2\bm{H}_{\T_{k+1},\T_{k+1}}(\bm{z}^k)\bm{h}^{k}_{\T_{k+1}} - \bm{J}_{\T_{k+1}}(\bm{z}^k)^{\top} \bm{J}_{\T_{k+1}}(\bm{h}^{k})\bm{h}^{k}_{\T_{k+1}}.
\end{align*}
According to the definition of $\bm{v}^{k}$, a simple calculation gives
\begin{align*} \label{vk-x}
\left\|\bm{v}^{k}_ {\T_{k+1}}-\bm{x}^\dag_{\T_{k+1}}\right\|_2 
~=&~ \left\|\bm{z}^{k}_ {\T_{k+1}}-\bm{x}^\dag_{\T_{k+1}}  - \mu^k \nabla{f}_{\T_{k+1}}(\bm{z}^k)\right\|_2 \cr
~\leq &~\underbrace{ \left\| \big(\bm{I} - 2\mu^k \bm{H}_{\T_{k+1},\T_{k+1}}(\bm{z}^k)\big) \bm{h}^{k}_{\T_{k+1}} \right\|_2}_{I_1} \cr
&~+\underbrace {\mu^k \left\|\bm{J}_{\T_{k+1}}(\bm{z}^k)^{\top} \bm{J}_{\T_{k+1}}(\bm{h}^{k})\bm{h}^{k}_{\T_{k+1}} \right\|_2 }_{I_2} .
\end{align*}
Next we will estimate ${I_1}$ and ${I_2}$ sequentially.
Since $\| \bm{h}^{k}\|_2 = \| \bm{z}^k - \bm{x}^\dag\|_2 \leq \gamma\norm{\bm{x}^\dag} $, we obtain
\[ 
(1-\gamma)\norm{\bm{x}^\dag} ~\leq~ \norm{\bm{z}^k} ~\leq~ (1+\gamma)\norm{\bm{x}^\dag}.
\]

For $I_1$: Denote $\lambda_1 := \lambda_{\max}\left( \bm{H}_{\T_{k+1},\T_{k+1}}(\bm{z}^k) \right)$, $\lambda_n := \lambda_{\min}\left( \bm{H}_{\T_{k+1},\T_{k+1}}(\bm{z}^k) \right)$. By \Cref{PHI} with $\bm{z} = \bm{z}^k$ and $\T = \T_{k+1}$, we have 
\begin{align}\label{G(z_k)}
 (2-\delta)(1-\gamma)^2\norm{\bm{x}^\dag}^2 &\leq \lambda_1 
 \leq (2+\delta)(1+\gamma)^2\norm{\bm{x}^\dag}^2 ,  \cr
 (\frac{1}{2}-\delta)(1-\gamma)^2\norm{\bm{x}^\dag}^2 &\leq \lambda_n \leq (\frac{1}{2}+\delta)(1+\gamma)^2\norm{\bm{x}^\dag}^2   .
\end{align}
Let $\mu^k \in \big(\frac{1}{(4-2\delta)(1-\gamma)^2\norm{\bm{x}^\dag}^2}, \frac{2}{(5+4\delta)(1+\gamma)^2\norm{\bm{x}^\dag}^2}\big)$, such that
\begin{align*}
0 ~\leq ~ 2\mu^k  \lambda_1 - 1  ~\leq~  1- 2\mu^k  \lambda_n .
\end{align*}
Consequently,
\begin{align*}
   ~&\big\| \bm{I} - 2\mu^k \bm{H}_{\T_{k+1},\T_{k+1}}(\bm{z}^k)\big\|_2 \cr
    \leq~& \max \{ |1- 2\mu^k  \lambda_n|, |1- 2\mu^k  \lambda_1| \}\cr
    \leq~& 1 - \mu^k(1 - 2\delta)(1-\gamma)^2\norm{\bm{x}^\dag}^2 ,
\end{align*}
let $l_0 = 1 - \mu^k(1 - 2\delta)(1-\gamma)^2 \norm{\bm{x}^\dag}^2$. Thus
\[
{I_1} ~\leq~   l_0 \norm{\bm{z}^k_{\T_{k+1}}  - \bm{x}^\dag_{\T_{k+1}} } .
\]

For $I_2$: An application of \Cref{PHI} with $\bm{z} = \bm{h}^{k}$ and $\T = \T_{k+1}$ implies
\[
\norm{\bm{J}_{\T_{k+1}}(\bm{h}^{k})} ~\leq~ \gamma\sqrt{2+\delta}\norm{\bm{x}^\dag},
\]
 together with \eqref{G(z_k)} we obtain
\[
{I_2} ~\leq~ (2+\delta)\gamma(1+\gamma)\mu^k\norm{\bm{x}^\dag}^2 \norm{\bm{z}^k_{\T_{k+1}} - \bm{x}^\dag_{\T_{k+1}} } .
\]	
Combining all pieces together, we have
\begin{align*}
\norm{\bm{v}_{\T_{k+1}}^{k}-\bm{x}^\dag_{\T_{k+1}}} 
~\leq~  I_1 +I_2 ~\leq~ \rho_0 \norm{\bm{z}_{\T_{k+1}}^k - \bm{x}^\dag_{\T_{k+1}}},
\end{align*}
where $\rho_0 = 1-\mu^{k}\norm{\bm{x}^\dag}^2(1 -2\delta  - (4 - 3\delta)\gamma - (1 + 3\delta)\gamma^2) $.
By \eqref{u-x-tk+1} we then have
\[
\norm{\bm{u}^{k}-\bm{x}^\dag}
~\leq~ \rho \norm{\bm{z}^{k}-\bm{x}^\dag},
\]
which completes the proof.
\end{IEEEproof}

Now we proceed to give the proof of \Cref{local convergence}, and we only consider the case where $\|\bm{z}^0 - \bm{x}^\dag\|_2 \leq \|\bm{z}^0 + \bm{x}^\dag\|_2$ since the case where $\|\bm{z}^0 + \bm{x}^\dag\|_2 \leq \|\bm{z}^0 - \bm{x}^\dag\|_2$ can be proved in a similar way.

\subsection{Proof of \Cref{local convergence}}\label{proof3}
We now give the proof of the first phase.
\begin{IEEEproof}
Given a vector $\bm{z}^{k}$ satisfying $\bm{z}^{k}\in\mathcal{E}(\gamma)$. $\{ \bm{z}^{k,l}\}_{l=0}^{L_0-1}$ are the iteration points produced by \Cref{algorithm} with $\bm{z}^{k,0} = \bm{u}^{k} $ and $\operatorname{supp}(\bm{z}^{k,l}) =  \cS_{k+1}$. The proof uses mathematical induction under the event \eqref{Phi-EPhi}. 
 Denote $\bm{h}^{k,l} := \bm{z}^{k,l} - \bm{x}^{\dag}$, $\cS^{\dag}: = \supp(\bm{x}^\dag)$, $\cS_k: = \supp(\bm{z}^k)$ and $\T_{k+1} := \cS_{k} \bigcup \cS_{k+1} \bigcup \cS^{\dag}$, which implies $| \T_{k+1}|\leq 3s$.  
According to \Cref{uk}, it holds that
\begin{align*}
    \left\|\bm{z}^{k,0}-\bm{x}^\dag\right\|_2 = \left\|\bm{u}^{k}-\bm{x}^\dag\right\|_2 \leq
 \rho\left\|\bm{z}^{k}-\bm{x}^\dag\right\|_2,
 \end{align*}
 where $\rho = \sqrt{\frac{3+\sqrt{5}}{2}}\big( 1-\mu^{k}\norm{\bm{x}^\dag}^2(1 -2\delta  - (4 - 3\delta)\gamma - (1 + 3\delta)\gamma^2)\big)$ and 
 $\mu^k \in \big(\frac{1}{(4-2\delta)(1-\gamma)^2\norm{\bm{x}^\dag}^2}, \frac{2}{(5+4\delta)(1+\gamma)^2\norm{\bm{x}^\dag}^2}\big)$. We can choose appropriate parameters to make sure that $\rho \in (0,1)$. A feasible choice of the parameters will be specified later.
 
Suppose that $\|\bm{z}^{k,l} - \bm{x}^\dag\|_2 \leq \|\bm{z}^{k} - \bm{x}^\dag\|_2$. According to the definition, we have
\begin{align}\label{J'F(u)}
&~\bm{J}_{\cS_{k+1}}(\bm{z}^{k,l})^{\top}\bm{F}(\bm{z}^{k,l})  \\
= &~\frac{1}{2m}\sum_{j=1}^m \Big((\bm{a}_{jR}^{\top}\bm{z}^{k,l})^2 + (\bm{a}_{jI}^{\top}\bm{z}^{k,l})^2 - (\bm{a}_{jR}^{\top}\bm{x}^\dag)^2- 
(\bm{a}_{jI}^{\top}\bm{x}^\dag)^2 \Big)  \nonumber  \\
&~ \cdot(\bm{a}_{jR,\cS_{k+1}}\bm{a}_{jR}^{\top}\bm{z}^{k,l} + \bm{a}_{jI,\cS_{k+1}}\bm{a}_{jI}^{\top}\bm{z}^{k,l}) \nonumber \\
= &~\bm{H}_{\cS_{k+1},\T_{k+1}}(\bm{z}^{k,l})\bm{h}^{k,l}_{\T_{k+1}} - \frac{1}{2}\bm{J}_{\cS_{k+1}}(\bm{z}^{k,l})^{\top}\bm{J}_{\T_{k+1}}(\bm{h}^{k,l})\bm{h}^{k,l}_{\T_{k+1}}. \nonumber 
\end{align}  
Since $\norm{\bm{h}^{k,l}} \leq \norm{ \bm{z}^{k} - \bm{x}^\dag } \leq \gamma \norm{\bm{x}^\dag} $, we have
\[    
(1-\gamma)\norm{\bm{x}^\dag} ~\leq~ \norm{ \bm{z}^{k,l}} ~\leq~ (1+\gamma)\norm{\bm{x}^\dag}.
\] 
By \Cref{PHI} with $\bm{z} = \bm{z}^{k,l}$ and $\T = \T_{k+1}$, we obtain
\begin{align}\label{1}
 \big\| \big(\bm{H}_{\cS_{k+1},\cS_{k+1}}(\bm{z}^{k,l} )\big )^{-1}\big\|_2 
~\leq~ \frac{2}{(1-2\delta)(1 - \gamma)^2\norm{\bm{x}^\dag}^2} ,
\end{align}
and
\begin{align}\label{2}
&~\big\|\big(\bm{H}_{\cS_{k+1},\cS_{k+1}}(\bm{z}^{k,l} )\big)^{-1} \bm{J}_{\cS_{k+1}}(\bm{z}^{k,l})^{\top}\big\|_2 \cr
\leq &~\frac{\sqrt{2}}{\sqrt{1-2\delta}(1 - \gamma)\norm{\bm{x}^\dag}}.
\end{align}
By \Cref{S-Sc} with $\bm{z} = \bm{z}^{k,l}$ and $\T = \T_{k+1}$ we have
\begin{align}\label{3}
\left\|\bm{H}_{\cS_{k+1}, \T_{k+1} \backslash \cS_{k+1}}(\bm{z}^{k,l} )\right\|_2 ~\leq~ \delta (1 + \gamma)^2\norm{\bm{x}^\dag }^2.
\end{align}
According to \Cref{PHI} with $\bm{z} = \bm{h}^{k,l} $ and $\T = \T_{k+1}$, then
\begin{align}\label{4}
\left\|\bm{J}_{\T_{k+1}}(\bm{h}^{k,l}) \right\|_2  ~\leq~ \sqrt{2 + \delta} \gamma\norm{\bm{x}^\dag}. 
\end{align}
We now decompose the term $\norm{\bm{z}^{k,l+1}-\bm{x}^\dag}$ into two parts and bound them respectively:
\begin{align*}
&~\left\|\bm{z}^{k,l+1}-\bm{x}^\dag\right\|_2^2 \cr
=&~ \left\| \bm{z}^{k,l+1}_{\cS^{\dag}\backslash\cS_{k+1}} - \bm{x}^\dag_{\cS^{\dag}\backslash\cS_{k+1}}\right\|_2^2  + \left\|\bm{z}^{k,l+1}_{\cS_{k+1}} - \bm{x}^\dag_{\cS_{k+1}}\right\|_2^2 \cr
=& ~ \left\| \bm{x}^\dag_{\cS^{\dag}\backslash\cS_{k+1}}\right\|_2^2 + \left\|\bm{z}^{k,l+1}_{\cS_{k+1}} - \bm{x}^\dag_{\cS_{k+1}}\right\|_2^2 .
\end{align*}   
Note that $\bm{x}^\dag_{\cS^{\dag}\backslash\cS_{k+1}}  $ is a sub-vector of $ \bm{z}^{k,0} - \bm{x}^\dag$, then we have
\begin{align}
\left\|\bm{x}^\dag_{\cS^{\dag} \backslash \cS_{k+1}}\right\|_2 ~\leq~\left\|\bm{z}^{k,0}-\bm{x}^\dag\right\|_2\leq
 \rho\left\|\bm{z}^{k}-\bm{x}^\dag\right\|_2.
\end{align}
According to update rule \eqref{iteration3} and \eqref{J'F(u)} we have	
\begin{align}  \label{[h_{k+1}]_S}
&~\bm{H}_{\cS_{k+1},\cS_{k+1}}(\bm{z}^{k,l} )(\bm{z}^{k,l+1}_{\cS_{k+1}}- \bm{x}^\dag_{\cS_{k+1}})  \cr
=&~\bm{H}_{\cS_{k+1},\cS_{k+1}}(\bm{z}^{k,l} )(\bm{z}^{k,l}_{\cS_{k+1}}- \bm{x}^\dag_{\cS_{k+1}}) - \bm{J}_{\cS_{k+1}}(\bm{z}^{k,l})^{\top}\bm{F}(\bm{z}^{k,l}) \cr
=&~\frac{1}{2}\bm{J}_{\cS_{k+1}}(\bm{z}^{k,l} )^{\top}\bm{J}_{\T_{k+1}}(\bm{h}^{k,l} )\bm{h}^{k,l}_{\T_{k+1}}  \cr
&~- \bm{H}_{\cS_{k+1}, \T_{k+1} \backslash \cS_{k+1}}(\bm{z}^{k,l} )\bm{h}^{k,l}_{\T_{k+1}\backslash\cS_{k+1}} ,
\end{align} 
thus 
\begin{align*}
&~\left\|\bm{z}^{k,l+1}_{\cS_{k+1}}- \bm{x}^\dag_{\cS_{k+1}}\right\|_2 \cr
=&~\left\|\big(\bm{H}_{\cS_{k+1},\cS_{k+1}}(\bm{z}^{k,l} )\big)^{-1} \bm{H}_{\cS_{k+1},\cS_{k+1}}(\bm{z}^{k,l} )(\bm{z}^{k,l+1}_{\cS_{k+1}}- \bm{x}^\dag_{\cS_{k+1}})  \right\|_2 \cr
\leq&~ \frac{1}{2}\left\|\big(\bm{H}_{\cS_{k+1},\cS_{k+1}}(\bm{z}^{k,l} )\big)^{-1}\bm{J}_{\cS_{k+1}}(\bm{z}^{k,l} )^{\top} \bm{J}_{\T_{k+1}}(\bm{h}^{k,l}) \bm{h}^{k,l}_{\T_{k+1}} \right\|_2  \cr
& + \left\| (\bm{H}_{\cS_{k+1},\cS_{k+1}}(\bm{z}^{k,l} ))^{-1}  \bm{H}_{\cS_{k+1}, \T_{k+1} \backslash \cS_{k+1}}(\bm{z}^{k,l} )\bm{h}^{k,l}_{\T_{k+1}\backslash\cS_{k+1}} \right\|_2 ,
\end{align*}
together with \eqref{1}, \eqref{2}, \eqref{3} and \eqref{4} yields that
\begin{align*}
&~\left\|\bm{z}^{k,l+1}_{\cS_{k+1}}- \bm{x}^\dag_{\cS_{k+1}}\right\|_2 \cr
\leq&~ \frac{\sqrt{2+\delta}\gamma}{\sqrt{2 - 4\delta}(1-\gamma)}\left\|\bm{z}^{k,l}_{\T_{k+1}} - \bm{x}^\dag_{\T_{k+1}} \right\|_2 \cr
 &~+ \frac{2\delta(1+\gamma)^2}{\left(1 - 2\delta\right)(1-\gamma)^2}\left\|\bm{z}^{k,l}_{\T_{k+1}\backslash\cS_{k+1}} - \bm{x}^\dag_{\T_{k+1}\backslash\cS_{k+1}}  \right\|_2  \cr
\leq&~ \alpha_0 \left\|\bm{z}^{k,l} - \bm{x}^\dag \right\|_2,
\end{align*}
where $\alpha_0 = \frac{\sqrt{(2 + \delta)(2 - 4\delta)}\gamma(1-\gamma) + 4\delta(1+\gamma)^2}{(2 - 4\delta)(1-\gamma)^2}$ . 
Combining all the terms, we obtain 
\begin{align*}
\norm{\bm{z}^{k,l+1}-\bm{x}^\dag}  &=\sqrt{\big\|\bm{x}^\dag_{\cS^\dag \backslash\cS_{k+1}}\big\|_2^2+\big\|\bm{z}^{k,l+1}_{\cS_{k+1}}-\bm{x}^\dag_{\cS_{k+1}}\big\|_2^2} \cr
&\leq \sqrt{\alpha_0^2\big\|\bm{z}^{k,l}-\bm{x}^\dag\big\|_2^2 +   \big\|\bm{z}^{k,0}-\bm{x}^\dag\big\|_2^2} \cr
&\leq \sqrt{\alpha_0^2 + \rho^2}\big\|\bm{z}^{k}-\bm{x}^\dag\big\|_2 .  
\end{align*}
Let $\alpha := \sqrt{\alpha_0^2 + \rho^2}$. We can choose parameters to make $\alpha, \rho \in (0,1)$. For example, we can set $\delta = 0.001$ and $\mu^k = 0.394/\|\bm{x}^\dag\|_2^2$, such that 
$\alpha, \rho \in (0,1)$ if provided $\gamma \leq 0.006$, which yields that $\norm{\bm{z}^{k,l}-\bm{x}^\dag} \leq \|\bm{z}^k-\bm{x}^\dag\|_2 $ hold for all $l = 0, 1,\cdots, L_0-1$. Therefore, we conclude $\norm{\bm{z}^{k+1}-\bm{x}^\dag} = \norm{\bm{z}^{k,L_0-1}-\bm{x}^\dag} \leq \alpha \left\|\bm{z}^k-\bm{x}^\dag\right\|_2 $.
\end{IEEEproof}

Next we give the proof of the second phase.
\begin{IEEEproof}
Given a vector $\bm{z}^{k}$ that satisfies $\bm{z}^{k}\in\mathcal{E}(\gamma)\cap\mathcal{E}(\frac{x^\dag_{\min}}{\norm{\bm{x}^\dag }})$ and $\{ \bm{z}^{k,l}\}_{l=0}^{L_0-1}$ are the iteration points produced by \Cref{algorithm} with $\bm{z}^{k,0} = \bm{u}^{k} $ and $\operatorname{supp}(\bm{z}^{k,l}) =  \cS_{k+1}$. The proof is based on mathematical induction under the event \eqref{Phi-EPhi}. 
 Let $\cS^{\dag}: = \supp(\bm{x}^\dag)$, $\cS_{k}: = \supp(\bm{z}^k)$ and $\bm{h}^{k,l} := \bm{z}^{k,l} - \bm{x}^\dag$. According to \Cref{uk}, it holds that
\begin{align*}
    \left\|\bm{z}^{k,0}-\bm{x}^\dag\right\|_2 = \left\|\bm{u}^{k}-\bm{x}^\dag\right\|_2 \leq
 \rho\left\|\bm{z}^{k}-\bm{x}^\dag\right\|_2 \leq \rho x^\dag_{\min},
 \end{align*}
 where $\rho = \sqrt{\frac{3+\sqrt{5}}{2}}\big( 1-\mu^{k}\norm{\bm{x}^\dag}^2(1 -2\delta  - (4 - 3\delta)\gamma - (1 + 3\delta)\gamma^2)\big)$ and  
$\mu^k \in \big(\frac{1}{(4-2\delta)(1-\gamma)^2\norm{\bm{x}^\dag}^2}, \frac{2}{(5+4\delta)(1+\gamma)^2\norm{\bm{x}^\dag}^2}\big)$. We can choose the appropriate parameters to ensure that $\rho \in (0,1)$. A feasible choice will be specified later.

 We confirm that $\cS^{\dag} \subseteq \cS_{k+1}$ holds when $\|\bm{z}^{k,0} - \bm{x}^\dag \|_2 < x^\dag_{\min} $, otherwise there would exist an index $j \in \cS^{\dag} \backslash \cS_{k+1} \neq \emptyset$,  such that $\|\bm{z}^{k,0}-\bm{x}^\dag\|_2 \geq|x_j| \geq x_{\min}^\dag$, which contradicts our assumption. 
 
 Suppose that $\|\bm{z}^{k,l} - \bm{x}^\dag\|_2 < \|\bm{z}^{k} - \bm{x}^\dag\|_2$, 
 then $\|\bm{z}^{k,l} - \bm{x}^\dag \|_2  <  \|\bm{z}^{k} - \bm{x}^\dag \|_2 <  \min \{ \gamma \|\bm x^\dag \|_2, x^\dag_{\min}\}$, which yields
\[ 
(1 - \gamma)\norm{\bm{x}^\dag} ~\leq~ \|\bm{z}^{k,l}\|_2 ~\leq~ (1 + \gamma)\norm{\bm{x}^\dag} .
\]

In terms of \Cref{PHI} with $\bm{z} = \bm{z}^{k,l}$ and $\T = \cS_{k+1}$, we obtain
\begin{align}\label{11}
&~\left\|(\bm{H}_{\cS_{k+1},\cS_{k+1}}(\bm{z}^{k,l} ))^{-1}\bm{J}_{\cS_{k+1}}(\bm{z}^{k,l} )^{\top}\right\|_2  \cr
\leq &~\frac{\sqrt{2}}{\sqrt{1-2\delta}(1 - \gamma)\norm{\bm{x}^\dag}}.
\end{align} 
\Cref{PHI} with $\bm{z} = \bm{h}^{k,l} $ and $\T = \cS_{k+1}$ implies
\begin{align}\label{13}
 \left\|\bm{J}_{\cS_{k+1}}(\bm{h}^{k,l}) \right\|_2  ~\leq~ \sqrt{2 + \delta
 } \left\|\bm{z}^{k,l} - \bm{x}^\dag \right\|_2. 
\end{align}
Similar to that in \eqref{[h_{k+1}]_S}, we can infer following equation
\begin{align*}
&~\bm{H}_{\cS_{k+1},\cS_{k+1}}(\bm{z}^{k,l} )(\bm{z}^{k,l+1}_{\cS_{k+1}}- \bm{x}^\dag_{\cS_{k+1}})   \cr
=&~\frac{1}{2}\bm{J}_{\cS_{k+1}}(\bm{z}^{k,l} )^{\top}\bm{J}_{\cS_{k+1}}(\bm{h}^{k,l} ) \bm{h}^{k,l}_{\cS_{k+1}},
\end{align*} 
together with \eqref{11} and \eqref{13} yields that
\begin{align*}
&~\big\|\bm{z}^{k,l+1}-\bm{x}^\dag \big\|_2  \cr
= &~\big\|\big(\bm{H}_{\cS_{k+1},\cS_{k+1}}(\bm{z}^{k,l} )\big)^{-1} \bm{H}_{\cS_{k+1},\cS_{k+1}}(\bm{z}^{k,l} ) (\bm{z}^{k,l+1}_{\cS_{k+1}}- \bm{x}^\dag_{\cS_{k+1}}\big)  \big\|_2 \cr
\leq&~ \frac{1}{2}\big\|\big(\bm{H}_{\cS_{k+1},\cS_{k+1}}(\bm{z}^{k,l} )\big)^{-1}\bm{J}_{\cS_{k+1}}(\bm{z}^{k,l} )^{\top}\bm{J}_{\cS_{k+1}}(\bm{h}^{k,l} ) \bm{h}^{k,l}_{\cS_{k+1}}\big\|_2  \cr
\leq &~ \frac{\sqrt{2+\delta}}{\sqrt{ 2 - 4\delta}(1- \gamma)\norm{\bm{x}^\dag }} \|\bm{z}^{k,l} - \bm{x}^\dag \|_2^2 \cr 
< &~ \beta\|\bm{z}^{k} - \bm{x}^\dag \|_2^2,
\end{align*}
where $\beta = \frac{\sqrt{2+\delta}}{\sqrt{ 2 - 4\delta}(1- \gamma)\norm{\bm{x}^\dag }}$. We can choose parameters $\delta$, $\mu^k$ and $\gamma$ to ensure that $\rho \in (0,1)$ and $\beta\cdot \gamma\| \bm{x}^\dag \|_2 \leq 1$. For example, we can set $\delta=0.001$ and $\mu^k = 0.394/\|\bm{x}^\dag\|_2^2$, such that $\rho \in (0,1)$ and $\beta\cdot \gamma\| \bm{x}^\dag \|_2 \leq 1$ if provided $\gamma \leq 0.006$, and therefore $\|\bm{z}^{k,l} - \bm{x}^\dag \|_2 < \|\bm{z}^{k} - \bm{x}^\dag \|_2$ hold for all $l = 0, 1, \cdots, L_0-1$. Then by induction, we have
\begin{align*}
  \big\|\bm{z}^{k+1}-\bm{x}^\dag \big\|_2 = \big\|\bm{z}^{k,L_0-1}-\bm{x}^\dag \big\|_2 < \beta \big\|\bm{z}^{k}-\bm{x}^\dag \big\|_2^2.
\end{align*}
\end{IEEEproof}

\subsection{Proof of Theorem~\ref{noisy case}}\label{proof5}
\begin{IEEEproof}
In noisy case, $y_j := y_j^{(\eta)}=|\langle\bm{a}_j, \bm{x}^\dag\rangle|^2 +\eta_j, j=1,\cdots,m$ with $\norm{\bm{\eta}} \leq C_\eta\|\bm{x}^\dag\|_2$. Given $\bm{z}^{k}$ satisfying $\bm{z}^{k} \in \mathcal{E}
(\gamma )$ and $\{ \bm z^{k,l} \}_{l=0}^{L_0-1}$ are the iteration points produced by \Cref{algorithm} with $\bm z^{k,0} = \bm u^k $ and $\supp (\bm z^{k,l}) = \mathcal{S}_{k+1}$. We use mathematical induction to prove this result under the event \eqref{Phi-EPhi}. Let $\cS^{\dag}: = \supp(\bm{x}^\dag)$, $\cS_k: = \supp(\bm{z}^k)$ and $\T_{k+1} := \cS_{k} \bigcup \cS_{k+1} \bigcup \cS^{\dag}$, which indicates that $|\T_{k+1} |\leq 3s$. 
Since
\begin{align*}
\bm{u}^{k}&=~ \mathcal{H}_s\big(\bm{z}^{k}- \mu^k\nabla f(\bm{z}^{k}) \big) \cr
&=~ \mathcal{H}_s\big(\bm{z}^{k}- \mu^k\bm{J}(\bm{z}^k)^{\top}\big(\bm{J}(\bm{z}^k) + \bm{J}(\bm{x}^\dag)\big)(\bm{z}^k - \bm{x}^\dag) \cr
&\qquad \qquad  - \mu^k\bm{J}(\bm{z}^k)^{\top}\bm{\eta} \big),
\end{align*}
using the same argument as the proof of the inequality in \Cref{uk}, we obtain
\begin{align} \label{u-v noise}
\norm{\bm{z}^{k,0}-\bm{x}^\dag } = &~  \norm{\bm{u}^{k}-\bm{x}^\dag } \cr
\leq &~ \sqrt{\frac{3+\sqrt{5}}{2}}\big(\bm{z}^k_{\T_{k+1}} - \bm{x}^\dag_{\T_{k+1}} -  \mu^k\nabla f_{\T_{k+1}}(\bm{z}^{k}) \big)  \cr
\leq&~ \rho \norm{\bm{z}^k - \bm{x}^\dag} + \zeta\norm{\bm{\eta}},
\end{align}
 where $\rho = \sqrt{\frac{3+\sqrt{5}}{2}}\big( 1-\mu^{k}\|\bm{x}^\dag\|_2^2(1 -2\delta  - (4 - 3\delta)\gamma - (1 + 3\delta)\gamma^2)\big)$, $\zeta =  \sqrt{\frac{3+\sqrt{5}}{2}}\mu^k\|\bm{x}^\dag\|_2\sqrt{2 + \delta}(1+\gamma)$ with $\mu^k \in \big(\frac{1}{(4-2\delta)(1-\gamma)^2\norm{\bm{x}^\dag}^2}, \frac{2}{(5+4\delta)(1+\gamma)^2\norm{\bm{x}^\dag}^2}\big)$. We can choose appropriate parameters such that $\rho \leq p\in(0,1)$ and $\zeta \leq q\in(0,1.62)$ hold for $l = 0$. A feasible choice will be specified later.

 Suppose that $\norm{\bm{z}^{k,l}-\bm{x}^\dag } \leq p\norm{\bm{z}^{k}-\bm{x}^\dag } + q\|\bm{\eta} \|_2$, since $\bm{z}^{k} \in \mathcal{E}(\gamma )$ and $\norm{\bm{\eta}} \leq C_{\eta}\norm{\bm{x}^\dag}$, then we have 
 \[
 \norm{\bm{z}^{k,l}-\bm{x}^\dag } \leq p \gamma \norm{\bm{x}^\dag} + C_{\eta} q \norm{\bm{x}^\dag} \leq \omega \norm{\bm{x}^\dag},
 \]
 where $\omega = p \gamma +  C_{\eta} q$, which implies that
 \[
 (1-\omega) \| \bm{x}^\dag  \|_2 \leq \| \bm{z}^{k,l}  \|_2 \leq (1+\omega) \| \bm{x}^\dag  \|_2 .
 \]
According to the update rule in \eqref{iteration3} and with the same argument as the proof of the first phase of \Cref{local convergence}, 
\begin{align}\label{z-x noise}		&~\left\|\bm{z}^{k,l+1}_{\cS_{k+1}}- \bm{x}^\dag_{\cS_{k+1}} \right\|_2 \cr
\leq&~ \bigg\|\bm{z}^{k,l}_{\cS_{k+1}}  -  \bm{x}^\dag_{\cS_{k+1}} -  \frac{1}{2}\big(\bm{H}_{\cS_{k+1}}(\bm{u}^{k} )\big)^{-1} \cr
&~\cdot  \bm{J}_{\cS_{k+1}}(\bm{z}^{k,l} )^{\top}\big(\bm{J}_{\T_{k+1}}(\bm{z}^{k,l} )\bm{z}^{k,l}_{\T_{k+1}} - \bm{J}_{\T_{k+1}}(\bm{x}^\dag )\bm{x}^\dag_{\T_{k+1}}\big)\bigg\|_2 \cr
&~ + \frac{1}{2}\norm{\bm{\eta} }\big\|(\bm{H}_{\cS_{k+1}}(\bm{z}^{k,l} ))^{-1}\bm{J}_{\cS_{k+1}}(\bm{z}^{k,l} )^{\top} \big\|_2 \cr
\leq&~p_0 \|\bm{z}^{k,l} - \bm{x}^\dag\|_2 + q_0\norm{\bm{\eta} },
\end{align}
where $p_0 = \frac{\sqrt{\left(2 + \delta\right)\left(2 - 4\delta\right)}\omega(1-\omega) + 4\delta(1+\omega)^2}{\left(2 - 4\delta\right)(1-\omega)^2} $ and $q_0 =\frac{1}{\sqrt{2-4\delta}(1-\omega)\norm{\bm{x}^\dag}}$.

Since $\bm{x}^\dag_{\cS^{\dag} \backslash \cS_{k+1}}$ is a sub-vector of $\bm{z}^{k,0}-\bm{x}^\dag$ and by $\sqrt{a^2+b^2} \leq|a|+|b|$, we obtain
\begin{align*}
\norm{\bm{z}^{k+1}-\bm{x}^\dag}
=&~\sqrt{\big\|\bm{x}^\dag_{\cS^{\dag} \backslash\cS_{k+1}}\big\|_2^2+\big\|\bm{z}^{k+1}_{\cS_{k+1}}-\bm{x}^\dag_{\cS_{k+1}}\big\|_2^2} \nonumber \cr
\leq&~ \big\|\bm{z}^{k,0}-\bm{x}^\dag\big\|_2  + p_0  \|\bm{z}^{k,l} - \bm{x}^\dag\|_2 + q_0 \norm{\bm{\eta} }  \cr
\leq&~ p_1\big\|\bm{z}^{k}-\bm{x}^\dag\big\|_2 + q_1\norm{\bm{\eta} } ,
\end{align*}
where $p_1 = p_0 p+ \rho$ and $q_1 = p_0q+ q_0 +\zeta$. 
There exist appropriate parameters $\delta, C_{\eta}, \mu^k$, such that $\rho, p_1\leq p\in (0,1)$ and $\zeta, q_1 \leq q\in(0,1.62)$, which implies $\left\|\bm{z}^{k,l}-\bm{x}^\dag\right\|_2 \leq p\left\|\bm{z}^{k}-\bm{x}^\dag\right\|_2 + q\norm{\bm{\eta} } $ hold for all $l = 0, 1, \cdots, L_0-1$. For example, one can choose $\delta = 0.0002$ , $C_{\eta} = 10^{-5}$ and step size $\mu^k = 0.394/\norm{\bm{x}^\dag}^2$ to make $\rho, p_1\leq p\in (0,1)$ and $\zeta, q_1 \leq q\in(0,1.62)$, if provided $\gamma \leq 0.001$. 
Then we have $\left\|\bm{z}^{k+1}-\bm{x}^\dag\right\|_2  = \left\|\bm{z}^{k, L_0-1}-\bm{x}^\dag\right\|_2 \leq p\left\|\bm{z}^{k}-\bm{x}^\dag\right\|_2 + q\norm{\bm{\eta} } $, which completes the proof.
\end{IEEEproof}

\subsection{Proof of \Cref{global}}\label{proof4}
 \begin{IEEEproof}
Define $\cS^{\dag}: = \supp(\bm{x}^\dag)$, $\cS_k: = \supp(\bm{z}^k)$ and $\T_{k+1} := \cS_{k} \bigcup \cS_{k+1} \bigcup \cS^{\dag}$, which implies $| \T_{k+1}|\leq 3s$.

According to \eqref{x^0}, with probability at least $1-8m^{-1}$, the $s$-sparse initial estimate $\bm{z}^0$ generated by initialization spectral method in \cite{jagatap2019sample}  satisfies $\bm{z}^{0}\in\mathcal{E}(\gamma)$.

Applying \Cref{local convergence} under the event \eqref{Phi-EPhi} throughout $K_{\max}$ iterations, then proceeding by induction, there exist an $\alpha \in (0,1)$ such that for any integer $0 \leq k < K_{\max}$, with probability at least $1-C_5 m^{-1}$, it holds that $\bm{z}^{k}\in\mathcal{E}(\gamma)$ and
\begin{align*}
\norm{\bm{z}^{k+1}-\bm{x}^{\dag}} ~\leq~ \alpha \left\|\bm{z}^{k}-\bm{x}^{\dag}\right\|_2.
\end{align*}

Let $K_0$ be the minimum integer such that
\begin{align*}
\gamma \alpha^{K_0}\norm{\bm{x}^\dag} ~<~ x^\dag_{\min} .
\end{align*}	
We can show that $\cS^{\dag} \subseteq \cS_{k}$ holds for all $k \geq K_0$ when $\|\bm{z}^{k}-\bm{x}^\dag\|_2 \leq \gamma \alpha^{k}\|\bm{x}^\dag\|_2\leq \gamma \alpha^{K_0}\|\bm{x}^\dag\|_2< x_{\min}^\dag$, otherwise there would exist an index $j \in \cS^{\dag} \backslash \cS_{k} \neq \emptyset$,  such that $\left\|\bm{z}^{k}-\bm{x}^\dag\right\|_2 \geq\left|x_j\right| \geq x_{\min}^\dag$, which contradicts with our assumption.  Thus, we have
\begin{align*}
K_0=\left\lfloor\frac{\log \big(\gamma \|\bm{x}^\dag\|_2/ x^\dag_{\min }\big)}{\log (\alpha^{-1})}\right\rfloor+1 \leq C_6 \log (\|\bm{x}^\dag\|_2 / x_{\min}^\dag)+C_7.
\end{align*}    
where $\lfloor\cdot\rfloor$ denotes the floor operation.
If $K_0 < K_{\max}$, then for $K_0 \leq k < K_{\max}$, according to the \Cref{local convergence} with the result that $\cS^{\dag} \subseteq \cS_{k+1}$, there exists $\beta$ such that $\beta \cdot \gamma\| \bm{x}^\dag \|_2 \leq 1$ and
\begin{align}\label{quadratic convergence rate}
\norm{\bm{z}^{k+1}-\bm{x}^\dag}  ~\leq~ \beta \left\|\bm{z}^{k} - \bm{x}^\dag\right\|_2^2. 
\end{align}

\end{IEEEproof}

\bibliographystyle{IEEEtran}  
\bibliography{GraHTP}

\vfill

\end{document}